\documentclass[11pt]{article}
\usepackage{amssymb,amsthm,hyperref}
\usepackage{amsmath}
\numberwithin{equation}{section}
\newtheorem{thm}{Theorem}[section]
\newtheorem{lem}{Lemma}[section]
\newtheorem{cor}{Corollary}[section]
\newtheorem{prop}{Proposition}[section]
\theoremstyle{definition}

\theoremstyle{remark}
\newtheorem{rem}{Remark}[section]
 \allowdisplaybreaks

\begin{document}
\title{A vacuum problem
for multidimensional compressible Navier-Stokes equations with
degenerate viscosity coefficients\thanks{This work is supported by
NSFC 10571158}}

\author{Ping Chen\thanks{Email: cp0804@hotmail.com} and Ting Zhang\thanks{Email: zhangting79@hotmail.com}
\\ Department of Mathematics, Zhejiang University
\\ Hangzhou 310027, PR China}

\date{}
\maketitle

\begin{abstract}
Local solutions of the multidimensional Navier-Stokes equations
for isentropic compressible flow are constructed with spherically
symmetric initial data between a solid core and a free boundary
connected to a surrounding vacuum state. The viscosity
coefficients $\lambda, \mu$ are proportional to $\rho^{\theta}$,
 $0<\theta<\gamma$, where $\rho$ is the density and $\gamma >
1$ is the physical constant of polytropic fluid. It is also proved
that no vacuum develops between the solid core and the free
boundary, and the free boundary expands with finite speed.
\end{abstract}
\section{Introduction}\label{sec1}
 We study compressible Navier-Stokes  equations
with degenerate viscosity coefficients in
$\mathbb{R}^{n}(n\geq2)$, which can be written in Eulerian
coordinates as
\begin{equation}\label{1.1}
\left \{
\begin{array}{l}
\partial_{t}\rho+\nabla\cdot(\rho \vec{u})=0, \\
\partial_{t}(\rho \vec{u})+\nabla\cdot(\rho \vec{u}\otimes \vec{u})+\nabla
P=\textrm{div}(\mu(\nabla \vec{u}+\nabla
\vec{u}^{\top}))+\nabla(\lambda \textrm{div}\vec{u})+\rho \vec{f},
\end{array} \right.
\end{equation}
where $\rho, P, \vec{u}=(u_{1}, \cdots, u_{n})$ and
$\vec{f}=(f_{1}, \cdots, f_{n})$ are the density,   pressure,
velocity and the given external force, respectively;
$\lambda=\lambda(\rho)$ and $\mu=\mu(\rho)$ are the viscosity
coefficients with the property $\lambda(0)=\mu(0)=0$.

In this paper,  we focus on the following initial-boundary value
problem for (\ref{1.1}).
 The initial conditions are
\begin{equation}\label{1.2}
(\rho, \vec{u})|_{t=0}=(\rho_{0}, \vec{u}_{0})(\textrm{x})
 =(\rho_{0}(r), u_{0}(r)\frac{\textrm{x}}{r}), 0<a\leq r\leq b,
\end{equation}
where $r=|\textrm{x}|=\sqrt{x_{1}^{2}+\cdots+x_{n}^{2}}$ and
$b>a>0$ are two constants, the boundary conditions are
\begin{eqnarray}
\vec{u}|_{r=a}=0, \label{1.3}\\
\rho|_{\textrm{x}\in\partial\Omega_{t}}=0,\label{1.4}
\end{eqnarray}
where $\partial\Omega_{t}=\psi(\partial\Omega_{0}, t)$ is a free
boundary. Here, $\partial\Omega_{0}=\{\textrm{x}\in\mathbb{R}^{n}
\ \big|\ |\textrm{x}|=b\}$ is the initial free boundary and $\psi$
is the flow of $u$:
\begin{equation*}
\left\{
\begin{array}{l}
  \partial_{t}\psi(\textrm{x},t)=\vec{u}(\psi(\textrm{x},t),t),\indent\textrm{x}\in\mathbb{R}^{n} \\
   \psi(\textrm{x},0)=\textrm{x}.
\end{array}
\right.
\end{equation*}

For the initial-boundary value problem (\ref{1.1})-(\ref{1.4})
with the spherically symmetric external force
$$\displaystyle \vec{f}(\textrm{x}, t)=f(r,
t)\frac{\textrm{x}}{r},$$ we are looking for a spherically
symmetric solution $(\rho, u)$:
$$\displaystyle\rho(\textrm{x}, t)
=\rho(r, t), \vec{u}(\textrm{x}, t) =u(r,
t)\frac{\textrm{x}}{r}.$$ Then $(\rho, u)$ is determined by
\begin{equation}\label{1.5}
\left \{
\begin{array}{lll}
\partial_{t}\rho+\partial_{r}(\rho u)+\frac{n-1}{r}\rho u&=&0, \\
\rho(\partial_{t}u+u\partial_{r}u)+\partial_{r}P
&=&(\lambda+2\mu)(\partial_{rr}^{2}u+\dfrac{n-1}{r}\partial_{r}u-\dfrac{n-1}{r^{2}}u)\\
&&+2\partial_{r}\mu\partial_{r}u+\partial_{r}
\lambda(\partial_{r}u+\dfrac{n-1}{r}u)+\rho f,
\end{array} \right.
\end{equation}
with the initial data
\begin{equation}\label{1.6}
(\rho, u)|_{t=0}=(\rho_{0}, u_{0})(r),\ 0<a\leq r\leq b,
\end{equation}
the fixed boundary condition
\begin{equation}\label{1.7}
u|_{r=a}=0,
\end{equation}
and the free boundary condition
\begin{equation}\label{1.8}
\rho(b(t), t)=0,
\end{equation}
where $b(0)=b$,  $b'(t)=u(b(t), t)$, $t>0$.

To handle this problem,  it is convenient to reduce the problem in
Eulerian coordinates $(r, t)$ to the problem in Lagrangian
coordinates $(x, t)$ moving with the fluid,  via the
transformation:
$$\displaystyle x=\int_{a}^{r}y^{n-1}\rho(y, t)dy, $$
then the fixed boundary $r=a$ and the free boundary $r=b(t)$
become
$$\tilde{a}=0\  \textrm{and}\  \tilde{b}(t)=\int_{a}^{b(t)}y^{n-1}\rho(y, t)dy
=\int_{a}^{b}y^{n-1}\rho_{0}(y)dy, $$ where
$\int_{a}^{b}y^{n-1}\rho_{0}(y)dy$ is the total mass initially.
Without loss of generality,  we can normalized it to 1,  so that
the region $\{(r, t)|a\leq r\leq b(t), 0\leq t<T\}$ under
consideration is transformed into the region $\{(x, t)|0\leq
x\leq1, 0\leq t<T\}$.

 Under the coordinate transformation,  the equations
(\ref{1.5})$\sim$(\ref{1.8}) are transformed into
\begin{equation}\label{1.9}
\left \{ \begin{array}{l}
\partial_{t}\rho=-\rho^{2}\partial_{x}(r^{n-1}u), \\
\partial_{t}u=r^{n-1}\partial_{x}\left(\rho(\lambda+2\mu)\partial_{x}(r^{n-1}u)-P\right)
              -2(n-1)r^{n-2}u\partial_{x}\mu+f(r, t),\\
r^{n}(x, t)=a^{n}+n\int_{0}^{x}\rho^{-1}(y, t)dy,
\end{array} \right.
\end{equation}
with the initial data
\begin{equation}\label{1.10}
(\rho, u)|_{t=0}=(\rho_{0}, u_{0})(x),  r|_{t=0}=r_{0}(x)
=\left(a^{n}+n\int_{0}^{x}\rho_{0}^{-1}(y)dy\right)^{\frac{1}{n}},
\end{equation}
and the boundary conditions:
\begin{eqnarray}
u(0, t)=0, \label{1.11}\\
\rho(1, t)=0.\label{1.12}
\end{eqnarray}

 It is a well-know fact that if we can solve the problem
(\ref{1.9})$\sim$(\ref{1.12}), then the free boundary problem
(\ref{1.5})$\sim$(\ref{1.8})
 has a solution.

 For simplicity of presentation,  we consider only the polytropic gas in this paper.
That is,  we assume
$$\left \{ \begin{array}{l}
P(\rho)=C\rho^{\gamma}, \gamma>1, \\
\mu(\rho)=c_{1}\rho^{\theta}, \\
\lambda(\rho)=c_{2}\rho^{\theta}, \theta>0. \\
\end{array} \right.$$
Without lose of generality, we assume $C=1$ . Additionally,  we
assume that the external force $ f(r, t)$ satisfies
\begin{equation}\label{1.14}
f(r, t)\in C^{1}([a, +\infty)\times[0, +\infty))\  \textrm{and}\
|f(r, t)|\leq \widetilde{f}(t), \ \forall\ r\geq a,\ t\geq0,
\end{equation}
with $\widetilde{f}\in C^{1}([0, \infty))$ and $\widetilde{f}\geq
0$.

The results in \cite{Hoff,Xin} show that the compressible
Navier-Stokes system with the constant viscosity coefficient has a
singularity at vacuum.  Considering the modified Navier-Stokes
system in which the viscosity coefficients depend on the density,
Liu, Xin and Yang in \cite{Liu} proved that such system is local
well-posed. It is motivated by the physical consideration that in
the derivation of the Navier-Stokes equations from the Boltzmann
equation through the Chapman-Enskog expansion to the second order,
cf. \cite{Grad}, the viscosity is a function of the temperature.
If we consider the case of isentropic gas, this dependence is
reduced to the dependence on the density function.

In this paper, we establish the existence, uniqueness and the
continuous dependence on initial data of the local solution to the
initial-boundary value problem (\ref{1.9})$\sim$(\ref{1.12}). The
important physical consequences of our results are that no vacuum
states can occur
 between the solid core and the free boundary,  the
interface separating the flow and vacuum propagates with finite
speed.

The main analytical difficulty to handle this problem is the
singularity of the solution near the free boundary. The most
difficult thing is to get the low bound of $\rho$, i.e.
$\rho\geq\frac{A}{3}(1-x)^{\alpha}$. Since $n\geq2$ and the
viscosity coefficient $\mu$ depends on $\rho$,
 the  nonlinear term  $r^{n-2}u\partial_{x}\mu$ in (\ref{1.9})$_2$ makes the analysis
  significantly different from  the one-dimensional
case \cite{FangZhang,Liu,Vong,Zhao,fang06}. For example, in the
proof of Proposition \ref{rho}, Lemmas \ref{ut2} and
\ref{ss-L3.11}, the nonlinear term do cause much trouble.  Another
difficult thing is to obtain the bound of $\|u\|_{L^{4m}}$. In
one-dimensional case, one can obtain the bound of $\|u\|_{L^{4m}}$
easily, but it's difficult in our case. We apply the inductive
method to prove it in Lemma \ref{u4m}. In the proof of Proposition
\ref{rho}, we must get $\int_x^1\frac{\partial_tu}{r^{n-1}}dy\leq
C\rho_0^{\theta}$. For this purpose, Yang-Zhao\cite{Zhao} needed
to obtain the $L^{4m}$ norm of $\partial_tu$ in one-dimensional
case. Here, we only need to get the $L^{2}$ norm of $\partial_tu$.
 In this paper, we use the
continuation method, the energy method and the pointwise estimate
techniques to deal with the singularity of the solution near the
free boundary, and obtain the bounded estimates of the solution.

Under some new \textit{a priori} estimates, using an effective
difference scheme to construct approximate solutions, we could
obtain the existence of the local solution with $0<\theta<\gamma$.
Using weighted function
$\rho_1^{l_1}\rho_2^{l_2}(\rho_1-\rho_2)^2$, $l_1+l_2=\theta-3$,
(not only $l_1=-1+\theta$ and $l_2=-2$ as in \cite{FangZhang1}),
we prove that the local solution is continuously dependent on the
initial data. To do this, we need $\alpha\theta\leq\frac{1}{2}$,
so $\theta$ can be any value in $(0, \gamma)$ when
$\alpha\leq\frac{1}{2\theta}$.

For related free boundary problems for isentropic fluids with
density-dependent viscosity ($\mu=\rho^\theta$), see
Fang-Zhang\cite{FangZhang}, Liu-Xin-Yang\cite{Liu},
Vong-Yang-Zhu\cite{Vong}, Yang -Zhao\cite{Zhao} and the references
cited therein. For the other one-dimensional or spherically
symmetric solutions of the Navier-Stokes equations for isentropic
compressible flow with a free boundary connected to a vacuum
state, see MatusuNecasova-Okada-Makino\cite{OM},
Nishida\cite{Nishida1,Nishida2}, Okada\cite{Okada},
Okada-Makino\cite{Okada2} and the references cited therein.  Also
see Lions\cite{Lions} and Vaigant\cite{Vaigant} for
multidimensional isentropic fluids. Especially,
Chen-Kratka\cite{Chen} studied the free boundary problem of
compressible heat-conducting flow with constant viscosity
coefficients and spherically symmetry initial data.

 Now we list some
assumptions on the initial data and  constants ($\gamma,\theta,
c_{1}, c_{2}$).
\par(A1)
Let $c_{1},  c_{2}, \theta,\gamma$ satisfy
$$c_1>0, 2c_1+nc_2>0, $$
\begin{equation}\label{1.17}
0<\theta<\gamma, \gamma>1.
\end{equation}
\par(A2)
There exist three positive constants $0<\alpha<1, A\geq B>0$
satisfying
\begin{eqnarray}\label{upalpha}
\frac{1}{2\gamma}<\alpha<\min\{\frac{3}{4\theta},
\frac{1}{\theta+1}\},\ \alpha(\theta-1)<\frac{1}{2},
\end{eqnarray}
and
\begin{equation}\label{1.15}
A(1-x)^{\alpha}\leq\rho_{0}(x)\leq B(1-x)^{\alpha}.
\end{equation}
\par(A3) Assume
    \begin{equation}\nonumber
    u_{0}\in  L^\infty([0, 1]),\ \rho_0^{1+\theta}(u_0)_x^2\in
    L^1([0,1]),
    \end{equation}
\begin{equation}\nonumber
(\rho_0^{\gamma})_x\in L^2([0,1]),\
(1-x)^{\alpha_{0}}(\rho_{0}^{\theta})_{x}^{2}\in L^{1}([0, 1]),
\end{equation}
\begin{equation}
\left((2c_1+c_2)\rho_{0}^{\theta+1}(r_{0}^{n-1}u_{0})_{x}\right)_{x}
-2c_1(n-1)\frac{u_0}{r_0}(\rho_0^{\theta})_x\in L^{2}([0,
1]),\label{1.16}
\end{equation}
where $1-2\alpha\theta<\alpha_{0}<\min\{1,
1+2\alpha-2\alpha\theta\}$, and $m$ is a integer satisfying
\begin{equation}\label{m}
m>\max\{\frac{1}{1+\alpha\theta-\alpha},
\frac{1}{2(1-\theta\alpha)}, \frac{1}{4-4\alpha}, 2\}.
\end{equation}

Under the above assumptions, our main results can be stated as
follows:
\begin{thm}[Existence]\label{thm} Under the assumptions (A1)$\sim$
(A3), there is a positive constant $T_{1}>0$ such that the free
boundary problem (\ref{1.9})$\sim$(\ref{1.12}) admits a weak
solution $(\rho, u, r)(x,t)$ on $[0, 1]\times[0, T_{1}]$ in the
sense that
$$\rho(x, t), u(x, t), r(x, t)\in L^\infty([0, 1]\times[0, T_{1}])
  \cap C^{1}([0, T_{1}];L^{2}([0, 1])),$$
        $$
        \frac{A}{3}(1-x)^{\alpha}\leq\rho(x, t)
     \leq3B(1-x)^{\alpha} ,\  (x,t)\in[0,1]\times[0,T_1],
       $$
$$\rho^{\theta+1}\partial_x(r^{n-1}u)\in L^{\infty}([0, 1]\times[0, T_{1}])\cap
C^{\frac{1}{2}}([0, T_{1}];L^{2}([0, 1])),$$
$$\partial_xr, \partial_xu \in L^{\infty}([0, T_{1}], L^{\lambda_0}([0, 1])),$$
where $\lambda_0>1$ is a constant which will be defined at
(\ref{zeta}), and the following equations hold:
 $$\partial_{t}\rho=-\rho^{2}\partial_{x}(r^{n-1}u),\indent \rho(x,0)=\rho_0, $$
 $$r^{n}(x, t)=a^{n}+n\int_{0}^{x}\rho^{-1}(y, t)dy,\ \partial_tr=u$$
 for almost all $x\in[0, 1]$,  any $t\in[0,T_1]$,
 \begin{eqnarray}\label{1.9_2}
\int_{0}^{T_{1}}\int_{0}^{1}u\phi_{t}dxdt
 &=&\int_{0}^{T_{1}}\int_{0}^{1}
\left(
  (2c_{1}+c_{2})\rho^{\theta+1}\partial_{x}(r^{n-1}u)-\rho^{\gamma}
 \right)
 \partial_{x}(r^{n-1}\phi)dxdt\nonumber\\
 &&-\int_{0}^{T_{1}}\int_{0}^{1}2c_{1}(n-1)\rho^{\theta}
                 \partial_{x}(\phi r^{n-2}u)dxdt\nonumber\\
&&-\int_{0}^{T_{1}}\int_{0}^{1}f\phi dxdt-\int_{0}^{1}u_{0}\phi(x,
0)dx
\end{eqnarray}
for any test function $\phi(x, t)\in C_{0}^{\infty}([0,
T_{1})\times(0, 1])$. Furthermore,  if
$\alpha\theta\leq\frac{1}{2}$, the solution satisfies
        \begin{equation}
          \rho\partial_x(r^{n-1}u)\in
          L^\infty([0,1]\times[0,T_1]),
        \end{equation}
    \begin{equation}
      (1-x)^{-\alpha}\rho\in C([0,T_1];L^\infty([0,1])).
    \end{equation}
\end{thm}

\begin{thm}[Continuous Dependence]\label{Continuous Dependence}
Under the assumptions of Theorem \ref{thm} and
$0<\alpha\theta\leq\frac{1}{2}$, if $(\rho_i, u_i, r_i)$ is a
solution to the system (\ref{1.9})$\sim$(\ref{1.12}) with the
initial data $(\rho_{0i}, u_{0i}, r_{0i})$, and satisfies
regularity conditions in Theorem \ref{thm}, $i=1,2$, then we have
\begin{eqnarray*}
&&\int^1_0\left((u_1-u_2)^2+\rho_1^{1-\theta}\rho_2^{2\theta-4}(\rho_1-\rho_2)^2+\rho_1^{\theta}\rho_2^{-1}(r_1-r_2)^2\right)dx\\
&\leq&Ce^{Ct}\int^1_0\left((u_{01}-u_{02})^2+\rho_{01}^{1-\theta}\rho_{02}^{2\theta-4}(\rho_{01}-\rho_{02})^2+\rho_{01}^{\theta}\rho_{02}^{-1}(r_{01}-r_{02})^2\right)dx,
\end{eqnarray*}
for all $0\leq t\leq T_1$.
\end{thm}
\begin{rem}
It is noticed that the set of initial data $(\rho_0, u_0, r_0)$
verifying all  assumptions in Theorem \ref{thm} and
\ref{Continuous Dependence} contains a quite general family of
functions. For example, if $\rho_0=C(1-x)^{\alpha}$ with the
exponent $\alpha$ satisfying (\ref{upalpha}), then it satisfies
all assumptions on the initial density.
\end{rem}

From  Theorem \ref{Continuous Dependence}, we can obtain the
uniqueness of the solution immediately.
\begin{thm}[Uniqueness]\label{uniq} Under the assumptions of Theorem
$\ref{Continuous Dependence}$, if the free boundary problem
(\ref{1.9})$\sim$(\ref{1.12}) has two weak solutions $(\rho_i,
u_i, r_i)$ ($i=1,2$) on $[0, 1]\times[0, T_{1}]$ as described in
Theorem $\ref{thm}$, then for all $t\in[0, T_{1}]$ and  almost all
$x\in(0,1)$, $(\rho_1,u_1,r_1)(x, t)=(\rho_2, u_2,r_2)(x, t)$.
\end{thm}
The rest of this paper is organized as follows. In Section
\ref{sec2}, we obtain some basic energy estimates,  the
 lower and upper bounds of the density. Then in Section
\ref{sec3}, we construct  approximate solutions by  the difference
scheme, and obtain the existence of the local solution. And in
Section \ref{sec4}, we prove Theorem \ref{Continuous Dependence}
by the energy method. At last, in Section \ref{sec5}, we prove
some useful lemmas, which are used in Section \ref{sec2}. In what
follows, we always use $C(C_i)$ to denote a generic positive
constant depending only on the initial data.

\section{Some \textit{a priori} estimates}\label{sec2}
In this section, for simplicity of presentation, we establish
certain \textit{a priori} estimates for smooth solutions to the
initial boundary value problem (\ref{1.9})$\sim$(\ref{1.12}).

\begin{prop}\label{energy}Let $(\rho, u, r)(x, t),\ x\in[0, 1],\ t\in[0, T]$, be a
solution of $(\ref{1.9}) \sim (\ref{1.12})$. Then under the
assumptions of Theorem $\ref{thm}$,  there exists a positive
constant $C_{1}=C_{1}
 (
  \|\rho_{0}\|_{L_{x}^{\infty}},
  \| u_{0}\|_{L_{x}^{2}}$,
  $\|\widetilde{f}\|_{L_t^\infty},  )$ such that,  for all  $t\in[0, T]$
\begin{eqnarray}\label{2.1}
&&\int_{0}^{1}(\frac{1}{2}u^{2}+\frac{1}{\gamma-1}\rho^{\gamma-1})(x,
t)dx
+\int_{0}^{t}\int_{0}^{1}\{(c_{2}+\frac{2}{n}c_{1})\rho^{\theta+1}
                 (\partial_{x}(r^{n-1}u))^{2}\nonumber\\
&&+\frac{2(n-1)}{n}c_{1}\rho^{\theta+1}(r^{n-1}\partial_{x}u
                -\frac{u}{r\rho})^{2}\}(x, s)dxds\nonumber\\
&\leq&
e^{t/2}\left(\int_{0}^{1}(\frac{1}{2}u_{0}^{2}+\frac{1}{\gamma-1}\rho_{0}^{\gamma-1})(x)dx
+\int_{0}^{t}\widetilde{f}^{2}(s)ds\right)\leq C_{1}e^{t/2}.
\end{eqnarray}
\end{prop}
\begin{proof} Multiplying $(\ref{1.9})_{1}$ by $\rho^{\gamma-2}$,  $(\ref{1.9})_{2}$ by $u$ and integrating it
over $[0, 1]\times[0, t]$,  using the boundary conditions
(\ref{1.11})$\sim$(\ref{1.12}) and Young's inequality, we have
\begin{eqnarray*}
&&\int_{0}^{1}(\frac{1}{2}u^{2}+\frac{1}{\gamma-1}\rho^{\gamma-1})(x,
t)dx
+\int_{0}^{t}\int_{0}^{1}\{(c_{2}+\frac{2}{n}c_{1})\rho^{\theta+1}
  (\partial_{x}(r^{n-1}u))^{2}\\
&&+\frac{2(n-1)}{n}c_{1}
  \rho^{\theta+1}(r^{n-1}\partial_{x}u-\frac{u}{r\rho})^{2}\}(x, s)dxds\\
&=&\int_{0}^{1}(\frac{1}{2}u_{0}^{2}+\frac{1}{\gamma-1}\rho_{0}^{\gamma-1})dx
   +\int_{0}^{t}\int_{0}^{1}fudxds\\
&\leq&\int_{0}^{1}(\frac{1}{2}u_{0}^{2}+\frac{1}{\gamma-1}\rho_{0}^{\gamma-1})dx
+\frac{1}{4}\int_{0}^{t}\int_{0}^{1}u^{2}dxds
+\int_{0}^{t}\widetilde{f}^{2}ds.
\end{eqnarray*}
Using Gronwall's inequality and (A1), we can obtain $(\ref{2.1})$
easily.
\end{proof}

\begin{lem}\label{r} Under the assumptions of Proposition  \ref{energy},  we have
for all $t\in[0,  T]$,
 \begin{eqnarray}
&&\quad\partial_{t}r(x, t)=u(x, t),\ \partial_{x}r(x,
t)=\frac{1}{r^{n-1}\rho}>0,\ r\geq a\label{2.2},
\\&&r^{\beta}\partial_{x}\rho^{\theta}(x, t)
=(r_{0}^{\beta}\partial_{x}\rho_{0}^{\theta})(x)
-\frac{\theta}{2c_{1}+c_{2}}((r^{1+\beta-n}u)(x,
t)-(r_{0}^{1+\beta-n}u_{0})(x))
\nonumber\\
    &&\quad+\frac{\theta}{2c_{1}+c_{2}}\int_{0}^{t}
   \left(
    -r^{\beta}\partial_{x}\rho^{\gamma}
    +r^{1+\beta-n}f
    +(1+\beta-n)u^{2}r^{\beta-n}\right)(x,s)ds\label{2.3},
\end{eqnarray}
where $\beta=\frac{2c_{1}\theta(n-1)}{2c_{1}+c_{2}}$ is a positive
constant.
\end{lem}
\begin{proof}From $(\ref{1.9})_3$, we can get (\ref{2.2}) immediately. From   $(\ref{1.9})_{1}$ and (\ref{2.2}),  we
have
\begin{eqnarray}
\partial_{t}(r^{\beta}\partial_{x}\rho^{\theta})
 &=&\beta r^{\beta-1}u\partial_{x}\rho^{\theta}
   +r^{\beta}\partial_{tx}\rho^{\theta},\label{2.5}\\
\partial_{t}\rho^{\theta}
  &=&-\theta\rho^{\theta+1}\partial_{x}(r^{n-1}u).\label{2.6}
\end{eqnarray}
From $(\ref{1.9})_{2}$ and $(\ref{2.5})\sim(\ref{2.6})$, we obtain
\begin{eqnarray*}
\partial_{t}(r^{\beta}\partial_{x}\rho^{\theta})
&=&\beta r^{\beta-1}u\partial_{x}\rho^{\theta}
-r^{\beta}\partial_{x}(\theta\rho^{\theta+1}\partial_{x}(r^{n-1}u))\\
&=&-\frac{\theta}{2c_{1}+c_{2}}(r^{1+\beta-n}(\partial_{t}u-f)
   +r^{\beta}\partial_{x}\rho^{\gamma}).
\end{eqnarray*}
Integrating both sides of the above equality in time $[0, t]$,
using integration by parts, we can obtain (\ref{2.3}) immediately.
\end{proof}

In Lemmas \ref{local-L2.3}-\ref{local-L3.4}, we use the
continuation method to estimate the lower bound of the density.
The key point is in the following proposition.

\begin{prop}\label{rho}Under the assumptions of Theorem $\ref{thm}$,
 there exists $T_{1}>0$, such that if
    \begin{equation}\label{rhobdry}
    \frac{1}{3}\rho_0(x)\leq \rho(x, t)\leq 3\rho_0(x),\ \forall\ t\in[0,
    T'],
    \end{equation}
where $T'\in(0, T_{1}]$ is any fixed positive constant, then we
have
    \begin{equation}\label{2.14}
    \frac{1}{2}\rho_0(x)\leq\rho(x, t)\leq 2\rho_0(x),\ \forall
    \ t\in[0, T'].
    \end{equation}
\end{prop}
\begin{rem}
Here $T_{1}$ can be defined by
\begin{equation}
T_{1}\equiv\min\{1, \overline{T}_{1}, \overline{T}_{2},
\overline{T}_{3}, \overline{T}_{4}\},\label{local-E3.8}
\end{equation}
where $\overline{T}_{1}$, $\overline{T}_{2}$, $\overline{T}_{3}$
and $\overline{T}_{4}$ are defined by (\ref{t2}), (\ref{t4}),
(\ref{2.40}) and (\ref{t1}), respectively. We assume that
$T_1\leq1$, since we only obtain the local existence of solution.

From (\ref{rhobdry}), we obtain
\begin{equation}\label{rup}
a\leq r\leq
\left(a^{n}+n\int_{0}^{x}(\frac{A}{3}(1-y)^{\alpha})^{-1}dy\right)^{\frac{1}{n}}
\leq R<\infty, \indent for\ all\  t\in[0, T'],
\end{equation}
where $R$ is a constant.
\end{rem}
Under the assumption (\ref{rhobdry}), we could obtain the
following lemmas. The proof of Lemmas \ref{u4m}$\sim$\ref{uxl1}
and Corollary \ref{cor} can be found in \textbf{Appendix 5}.
\begin{lem}\label{u4m}Under the assumptions of Theorem $\ref{thm}$ and (\ref{rhobdry}),
 there exists a positive constant $C_{2}=C_{2}
           (a,
            \|\rho_{0}\|_{L_{x}^{\infty}},
            \| u_{0}\|_{L_{x}^{4m}}, \
            \|\widetilde{f}\|_{L_{t}^{\infty}}, \
            \ C_{1}
           )$,
such that for all $t\in[0, T']$ and $k=1, 2, \cdots, 2m$,
\begin{eqnarray}\label{u4m1}
    \int_{0}^{1}u^{2k}(x, t)dx
   +\int_{0}^{t}\int_{0}^{1}
\rho^{\theta+1}u^{2k-2}r^{2n-2}(\partial_{x}u)^{2} dxds \leq
C_{2}.
\end{eqnarray}
\end{lem}
\begin{lem}\label{ul2} Under the assumptions of Lemma \ref{u4m},  there exists a
positive constant $C_{3}=C_{3}(a,
R,A,B,\|\rho_0^{1+\theta}(u_0)_x^2\|_{L^1}$, $
\|\widetilde{f}\|_{L_t^\infty},  C_{1})$ such that, for all
$t\in[0, T']$,
\begin{equation}\label{2.16}
\int_{0}^{1}|u(x, t)-u_{0}(x)|^{2}dx\leq C_{3}t.
\end{equation}
\end{lem}
\begin{cor}\label{cor}Under the assumptions of Lemma \ref{u4m},
there exists a positive constant $C_{4}=C_{4}(a, R$, $
\|u_0\|_{L_x^\infty},  C_{2}, C_{3})$ such that, for all $t\in[0,
T']$,
\begin{equation}\label{2.21}
\left|\int_{0}^{t}\int_{x}^{1}\frac{\partial_{t}u}{r^{n-1}}dxds\right|\leq
C_{4}t^{\frac{1}{2(2m-1)}}(1-x)^{\frac{2m-1}{2m}}.
\end{equation}
\end{cor}
\begin{lem}\label{weightedrho}Under the assumptions of Lemma \ref{u4m},
there exists a positive constant $C_{5}=C_{5}(a, R$, $
\|(1-x)^{\alpha_0}(\rho_0^{\theta})_x^2\|_{L^1},
\|\widetilde{f}\|_{L_t^\infty},  C_{2})$, such that, for $t\in[0,
T']$,
\begin{equation}
\int_{0}^{1}(1-x)^{\alpha_{0}}(\partial_{x}\rho^{\theta}(x,
t))^{2}dx\leq C_{5},\label{local-E3.13}
\end{equation}
where $\alpha_{0}>1-2\alpha\theta$.
\end{lem}

\begin{lem}\label{ut2}Under the assumptions of Lemma \ref{u4m}, there exist a positive
constant $C_{6}$ which is defined by ($\ref{C12}$) and a positive
constant $C_7$ such that, for all $t\in[0, T']$,
\begin{eqnarray}
&\displaystyle\int_{0}^{1}(\partial_{t}u)^{2}dx
+\int_0^t\int_0^1\left[\rho^{\theta+1}r^{2n-2}(\partial_{tx}u)^2
+\rho^{\theta-1}r^{-2}(\partial_tu)^2\right]dxds \leq
C_{6},\label{2.20}
&\\
&\displaystyle\int^1_0\left(\rho^{1+\theta}(\partial_xu)^2+\rho^{\theta-1}u^2\right)dx+
    \int^t_0\int^1_0(\partial_tu)^2dxds\leq C_{7}&\label{2.155}\\&
\displaystyle\int_{0}^{1}\rho^{\theta+3}(\partial_{x}u)^{4}dx\leq
C_{7}.&\label{2.15}
\end{eqnarray}
\end{lem}
\begin{lem}\label{uxl1}Under the assumptions of Lemma \ref{u4m},
there exists a positive constant $C_{8}$ such that, for all
$t\in[0, T']$,
\begin{eqnarray}
\| u\|_{L^{\infty}([0, 1]\times[0, T'])}\leq C_8,\label{b2}\\
 \int_{0}^{1}|\partial_{x}u|^{\lambda_0}dx\leq C_{8},\label{c15}
\end{eqnarray}
where $\lambda_0$ is a constant satisfying:
\begin{equation}\label{zeta}
1<\lambda_0<\min\{\frac{4m}{4m\alpha+1},\frac{1}{\alpha(1+\theta)}\}.
\end{equation}
\end{lem}
\begin{rem}
Since $m>\frac{1}{4-4\alpha}$ and $\alpha<\frac{1}{\theta+1}$, we
have $1<\min\{\frac{4m}{4m\alpha+1},\frac{1}{\alpha(1+\theta)}\}$
and the set of $\lambda_0$ is not empty.
\end{rem}

Now, we turn to prove \textbf{Proposition \ref{rho}}:
\begin{proof}From (\ref{1.9}), we have
\begin{eqnarray*}
    &&\rho^{\theta}(x,
    t)+\frac{\theta}{2c_{1}+c_{2}}\int_{0}^{t}\rho^{\gamma}(x,
    s)ds=\rho_{0}^{\theta}(x)
         +\frac{\theta}{2c_{1}+c_{2}}
          \int_{0}^{t}\int_{x}^{1}\left(\frac{\partial_{t}u}{r^{n-1}
          }\right)(y,s)dyds\\
        &&-\frac{2c_{1}(n-1)\theta}{2c_{1}+c_{2}}
    \int_{0}^{t}\left(\frac{u\rho^{\theta}}{r}\right)(x,s)ds
    -\frac{2c_{1}(n-1)\theta}{2c_{1}+c_{2}}
    \int_{0}^{t}\int_{x}^{1}\left(\rho^{\theta}\partial_{x}\left(
    \frac{u}{r}\right)\right)(y,s)dyds\\
    &&
    -\frac{\theta}{2c_{1}+c_{2}}\int_{0}^{t}\int_{x}^{1}(fr^{1-n})(y,
                s)dyds.
    \end{eqnarray*}
Using (\ref{2.21}), we get
\begin{eqnarray}\label{1}
&&\rho^{\theta}(x,
t)+\frac{\theta}{2c_{1}+c_{2}}\int_{0}^{t}\rho^{\gamma}(x, s)ds
\nonumber\\
&\geq& \rho_{0}^{\theta}(x)
-\frac{\theta}{2c_{1}+c_{2}}C_{4}t^{\frac{1}{2(2m-1)}}(1-x)^{\frac{2m-1}{2m}}-
\frac{2c_{1}(n-1)\theta}{2c_{1}+c_{2}}
\int_{0}^{t}\frac{u\rho^{\theta}}{r}ds
\nonumber\\
&&-\frac{2c_{1}(n-1)\theta}{2c_{1}+c_{2}}
\int_{0}^{t}\int_{x}^{1}\rho^{\theta}\partial_{x}\left(\frac{u}{r}\right)dyds
-\frac{\theta}{2c_{1}+c_{2}}a^{-n+1}\|\widetilde{f}\|_{L_t^\infty}(1-x)t\nonumber\\
&\geq& \rho_{0}^{\theta}(x) -C_{9}
t^{\frac{1}{2(2m-1)}}\rho_{0}^{\theta}-
\frac{2c_{1}(n-1)\theta}{2c_{1}+c_{2}}
\int_{0}^{t}\frac{u\rho^{\theta}}{r}ds
\nonumber\\
&&-\frac{2c_{1}(n-1)\theta}{2c_{1}+c_{2}}
\int_{0}^{t}\int_{x}^{1}\rho^{\theta}\partial_{x}\left(\frac{u}{r}\right)dyds
-C_{10}(1-x)t,
\end{eqnarray}
where $m\geq\frac{1}{2(1-\theta\alpha)}$ and
$C_{9}:=\frac{\theta}{2c_1+c_2}C_{4},
C_{10}:=\frac{\theta}{2c_1+c_2}a^{-n+1}\|\widetilde{f}\|_{L_t^\infty}.$

Using Lemma \ref{uxl1}, (\ref{2.2}), (\ref{rhobdry}) and
$0<\alpha<1$, we have
    \begin{eqnarray}\label{2}
   \left|\int_{0}^{t}\frac{u(x,s)\rho^{\theta}(x,s)}{r(x,
    s)}ds\right|
    \leq\frac{\left(3\rho_0\right)^{\theta}}{a}\int_{0}^{t}|u(x, s)|ds
        \leq C_{8}\frac{(3\rho_0)^{\theta}}{a}t
        :=C_{11}\rho_{0}^{\theta}t,
    \end{eqnarray}
\begin{eqnarray}\label{3}
\left|\int_{0}^{t}\int_{x}^{1}\rho^{\theta}\frac{\partial_{x}u}{r}dxds\right|
\leq\frac{\left(3\rho_0\right)^{\theta}}{a}\int_{0}^{t}\int_{x}^{1}|\partial_{x}u|dxds
\leq C_{8}\frac{(3\rho_0)^{\theta}}{a}t
:=C_{12}(\rho_{0})^{\theta}t,
\end{eqnarray}
and
    \begin{eqnarray}\label{4}
   \left|\int_{0}^{t}\int_{x}^{1}\rho^{\theta}\frac{u\partial_{x}r}{r^{2}}dxds\right|
     &=&\left|\int_{0}^{t}\int_{x}^{1}\rho^{\theta-1}r^{-n-1}udxds\right|\nonumber\\
            &\leq&a^{-n-1}(3\rho_0)^\theta
            \int_0^t\|u\|_{L^{\infty}}\|3(\rho_0)^{-1}\|_{L^1}ds\nonumber\\
                    &\leq& Ct(\rho_{0})^{\theta}\int^1_0(1-x)^{-\alpha} dx
                    \leq      C_{13}(\rho_{0})^{\theta}t.
    \end{eqnarray}
We define
\begin{equation}\label{t2}
\overline{T}_{1}=\left(\frac{(1
                  -(\frac{2}{3})^{\theta})A^\theta
                  }
  {C_{9}+C_{10}+C_{11}+C_{12}+C_{13}}\right)^{2(2m-1)}.
\end{equation}
From $(\ref{1})\sim(\ref{4})$, we have
\begin{equation}\label{2.36}
\rho^{\theta}(x,
t)+\frac{\theta}{2c_{1}+c_{2}}\int_{0}^{t}\rho^{\gamma}(x, s)ds
\geq(\frac{2}{3}\rho_0)^{\theta},
\end{equation}
for all $0\leq t\leq T'\leq T_{1}\leq \overline{T}_{1}.$

To get a lower bound on $\rho(x, t)$ from (\ref{2.36}), we need to
get a upper bound of the term $\int_{0}^{t}\rho^{\gamma}(x, s)ds$
for sufficiently small $t$.  For this purpose,  we set
\begin{equation}\label{2.37}
Z(t)=\int_{0}^{t}\rho^{\gamma}(x, s)ds.
\end{equation}
Similar to the proof of (\ref{2.36}), we can deduce that $Z(t)$
satisfies the following differential inequality from
$(\ref{1.9})$:
\begin{equation}\label{2.38}
\left(Z'(t)\right)^{\frac{\theta}{\gamma}}+\frac{\theta}{2c_{1}+c_{2}}
Z(t) \leq \rho_0^\theta+(C_{9}
          +C_{10}
          +C_{11}
          +C_{12}
          +C_{13})t^{\frac{1}{2(2m-1)}}(1-x)^{\theta\alpha},
\end{equation}
for $t\in[0,T']$. Define
\begin{equation}\label{t4}
\overline{T}_{2}=\left(\frac{(2^{\theta}-1)A^\theta}
                       {C_{9}+C_{10}+C_{11}+C_{12}+C_{13}}
                        \right)^{2(2m-1)}.
\end{equation}
For all $t\leq T'\leq \overline{T}_{2}$, we have
    $$
    \left(Z'(t)\right)^{\frac{\theta}{\gamma}}+\frac{\theta}{2c_{1}+c_{2}}
Z(t) \leq C_{14}(1-x)^{\alpha\theta}
    $$
and
\begin{equation*}
\rho(x,t)\leq2\rho_0(x),\ \forall\ x\in[0,1],\ t\in[0,T'].
\end{equation*}
This is the second inequality of (\ref{2.14}).

Noticing $Z(0)=0$ and $0<\theta<\gamma$, we finally deduce from
(\ref{2.38}) that
\begin{equation}\label{2.39}
-\frac{\theta}{2c_{1}+c_{2}} Z(t)\geq C_{14}(1-x)^{\theta\alpha}
       (
         (
          1+\frac{\gamma-\theta}{2c_{1}+c_{2}}
          (C_{14}(1-x)^{\theta\alpha})^{\frac{\gamma-\theta}{\theta}}t
         )^{\frac{\theta}{\theta-\gamma}}-1
       ).
\end{equation}
By choosing $\overline{T}_{3}>0$ sufficiently small such that
\begin{equation}\label{2.40}
C_{14}
       \left(
         1-\left(
          1+\frac{\gamma-\theta}{2c_{1}+c_{2}}
          (C_{14}(1-x)^{\theta\alpha})^{\frac{\gamma-\theta}{\theta}}\overline{T}_{3}
         \right)^{\frac{\theta}{\theta-\gamma}}
       \right)\leq\left((\frac{2}{3})^{\theta}-(\frac{1}{2})^{\theta}\right)A^{\theta},
\end{equation}
we can get from (\ref{2.39})$\sim$(\ref{2.40}) that
\begin{equation}\label{2.41}
-\frac{\theta}{2c_{1}+c_{2}}\int_{0}^{t}\rho^{\gamma}(x, s)ds
\geq-\left((\frac{2}{3})^{\theta}-(\frac{1}{2})^{\theta}\right)\rho_0^{\theta}(x).
\end{equation}
Inserting (\ref{2.41}) into (\ref{2.36}),  we can arrive at
\begin{equation}
\rho^{\theta}\geq(\frac{1}{2}\rho_0)^{\theta},
\end{equation}
provided that $0\leq t\leq T'\leq T_{1}$. This completes the proof
of Proposition \ref{rho}.
\end{proof}

Using the pointwise bound for $\rho$  (\ref{rhobdry}), we can also
derive other estimates for $(\rho,u)$ in the following
propositions.

\begin{prop}
Under the conditions in Lemma \ref{u4m}, we have
\begin{eqnarray}\label{rhoc}
     &  \displaystyle \int_{0}^{1}|\partial_t\rho(x, t)|^{2}dx \leq C,
   \\\label{xc}
      &    \displaystyle  \int^t_0\int_{0}^{1}\left|\partial_t\left(\rho^{\theta+1}\partial_{x}(r^{n-1}u)\right)(x,
            \eta)\right|^{2}dxd\eta
            \leq C,\ t\in[0,T'].
            \end{eqnarray}
\end{prop}
\begin{proof}From equation (\ref{1.9})$_1$ and H\"{o}lder's
inequality, we have
    \begin{eqnarray*}
   \int_{0}^{1}\left|\rho_{t}(x, t)\right|^{2}dx
        &=&\int_{0}^{1}\left|\left(\rho^{2}\partial_{x}(r^{n-1}u)\right)(x, t)\right|^{2}dx\\
                &\leq&C\int_{0}^{1}\left(
                  \rho^{2}r^{-2}u^2+\rho^{4}r^{2n-2}(\partial_{x}u)^{2}
                 \right)(x, t)dx.
    \end{eqnarray*}
Using Lemma \ref{ut2}, (\ref{2.1}), (\ref{rhobdry}), (\ref{rup}),
(\ref{2.15}) and $\alpha<\frac{1}{1+\theta}$, we have
  \begin{eqnarray*}
    &&\quad\int_{0}^{1}\rho^{2}r^{-2}u^2(x, t)dx \leq C\int_{0}^{1}u^{2}(x,
    t) dx \leq C,\\
        &&\quad\int_{0}^{1}\rho^{4}r^{2n-2}(\partial_{x}u)^{2}(x,
        t) dx\leq
        C\int_{0}^{1}(\rho^2\partial_{x}u)^{2}(x,
        t) dx\\
              &  &\leq C\left(\int_0^1\rho^{\theta+3}(\partial_xu)^4dx+\int_0^1\rho^{5-\theta}dx\right)
                \leq C.
        \end{eqnarray*}
From above, we get (\ref{rhoc}) immediately.

Now, we turn to prove (\ref{xc}):
\begin{eqnarray*}
  &&\int^t_0\int_{0}^{1}\left|\partial_{t}\left(\rho^{\theta+1}
               \partial_{x}(r^{n-1}u)\right)(x, t)\right|^{2}dxd\eta\\
  &\leq&C\int_{0}^{t}\int_{0}^{1}\left(\rho^{2\theta+4}
                    [\partial_{x}(r^{n-1}u)]^{4}+\rho^{2\theta+2}
                              \left(\partial_{tx}(r^{n-1}u)\right)^{2}\right)(x,
                              \eta)dxd\eta.
\end{eqnarray*}
Denote
\begin{eqnarray*}
D_{1}&=&\int_{0}^{t}\int_{0}^{1}
             \left(
               \rho^{2\theta+4}\left(\partial_{x}(r^{n-1}u)\right)^{4}
             \right)(x, \eta)dxd\eta, \\
D_{2}&=&\int_{0}^{t}\int_{0}^{1}
             \left(
             \rho^{2\theta+2}\left(\partial_{tx}(r^{n-1}u)\right)^{2}
             \right)(x, \eta)dxd\eta.
\end{eqnarray*}
From Lemmas \ref{u4m}, \ref{ut2}, (\ref{rhobdry}) and (\ref{rup}),
we have
\begin{eqnarray*}
D_{1}&\leq&C\int_{0}^{t}\int_{0}^{1}
                \rho^{2\theta+4}\left(\rho^{-4}r^{-4}u^{4}
                         +r^{4n-4}(\partial_{x}u)^{4}\right)dxd\eta\\
     &\leq&C\int_{0}^{t}\int_{0}^{1}
                \left(u^{4}+\rho^{\theta+3}(\partial_{x}u)^{4}\right)dxd\eta\\
     &\leq&C
\end{eqnarray*}
and
\begin{eqnarray*}
D_{2}&\leq&C\int_{0}^{t}\int_{0}^{1}
             (
              \rho^{2\theta+2}(\partial_{x}(r^{n-2}u^{2}))^{2}
              +\rho^{2\theta+2}(\partial_{x}(r^{n-1}\partial_{t}u))^{2}
             )dxd\eta\\
&\leq&C\int_{0}^{t}\int_{0}^{1}\left[\rho^{2\theta+2}r^{-4}\rho^{-2}u^{4}
             +\rho^{2\theta+2}r^{2n-4}(\partial_{x}u)^{2}u^{2}\right]dxd\eta
+CC_{6}\\
&\leq& C.
\end{eqnarray*}
This completes the proof.
\end{proof}
\begin{prop}\label{local-P3.4} Under the assumptions of
Lemma \ref{u4m},  we get
\begin{eqnarray}\label{b1}
&\displaystyle\|
    [\rho^{\theta+1}\partial_{x}(r^{n-1}u)](x, t)
\|_{L^{\infty}([0, 1]\times[0, T'])}\leq C,
\\&\label{b3}
\displaystyle\int_{0}^{1}\left|\partial_{x}\left(\rho^{\theta+1}\partial_{x}(r^{n-1}u)\right)\right|dx\leq
C,
\\\label{b4}
   &\displaystyle \int_{0}^{1}|\partial_{x}\rho|dx\leq C,
    \ \textrm{ for all}\ t\in[0,T'].
\end{eqnarray}
\end{prop}
\begin{proof}From $(\ref{1.9})_{2}$, we have
\begin{eqnarray*}
\partial_{x}[\rho^{\theta+1}\partial_{x}(r^{n-1}u)]
& =&\frac{\partial_{t}u}{(2c_{1}+c_{2})r^{n-1}}
  +\frac{\partial_{x}\rho^{\gamma}}{2c_1+c_2}
  \\
  &&+\frac{2c_{1}(n-1)u\partial_{x}\rho^{\theta}}{(2c_{1}+c_{2})r}
  -\frac{f}{(2c_{1}+c_{2})r^{n-1}}.
\end{eqnarray*}
From (\ref{1.14}) and Lemmas \ref{ut2}$\sim$\ref{uxl1}, using
similar arguments as that in (\ref{2})$\sim$(\ref{4}), we get
\begin{eqnarray*}
        &&\left|\rho^{\theta+1}\partial_{x}(r^{n-1}u)\right|\\
                    &\leq&\frac{1}{(2c_{1}+c_{2})a^{n-1}}\int_{x}^{1}\left|\partial_{t}u\right|dx
                      +\frac{1}{2c_1+c_2}\left|\rho^{\gamma}\right|\\
        &&+\frac{2c_{1}(n-1)}{(2c_{1}+c_{2})}\left(\rho^\theta\frac{|u|}{r}
        +\left|\int_{x}^{1}\rho^{\theta}\partial_{x}\left(
        \frac{u}{r}\right)dx\right|\right)
        +\frac{1}{(2c_{1}+c_{2})a^{n-1}}\int_{x}^{1}|f|dx\\
                & \leq& C.
    \end{eqnarray*}
 Similarly, using (\ref{local-E3.13}) and
$\alpha_0<1$, we have
        \begin{eqnarray*}
       \int_{0}^{1}\left|\partial_{x}\left(\rho^{\theta+1}\partial_{x}(r^{n-1}u)\right)\right|dx
                &\leq&C+C\int^1_0|\partial_x(\rho^\theta)|dx\\
                    &\leq&C+C\int^1_0(1-x)^{\alpha_0}
                    |(\rho^\theta)_x|^2dx\!+\!C\int^1_0(1\!-\!x)^{-\alpha_0}dx\\
        &\leq& C+CC_5+\int^1_0(1-x)^{-\alpha_0}dx\leq C.
        \end{eqnarray*}

From (\ref{rhobdry}), (\ref{local-E3.13}),
$\alpha_0<1+2\alpha-2\alpha\theta$ and H\"{o}lder's inequality, we
get
    \begin{eqnarray*}
      \int^1_0|\rho_x|dx&\leq& C \left(\int^1_0(1-x)^{\alpha_0}
      (\rho^\theta)_x^2dx
      \right)^\frac{1}{2}\left(\int^1_0(1-x)^{-\alpha_0}
      \rho^{2-2\theta}dx
            \right)^\frac{1}{2}\\
                    &\leq&C\left(\int^1_0(1-x)^{-{\alpha_0}+2\alpha-2\alpha\theta}dx
            \right)^\frac{1}{2}\leq C.
    \end{eqnarray*}
This completes the proof.
\end{proof}

\begin{lem}
 Under the assumptions of
Lemma \ref{u4m}, we have
    $$
    \|r(\cdot,t)\|_{L^\infty\cap W^{1,\lambda_0}}\leq C,
    \ \forall \ t\in[0,T'].
    $$
\end{lem}

\begin{lem}\label{lem}
 Under the assumptions of
Lemma \ref{u4m} and $0<\alpha\theta\leq\frac{1}{2}$, there exists
a constant $C_{15}>0$ such that
    \begin{equation}\label{lembdry}
    \|\rho \partial_xu\|_{L^{\infty}([0, 1]\times[0, T'])}\leq
    C_{15},
    \end{equation}
        \begin{equation}
          \|\rho_0^{-1}\rho_t\|_{L^\infty([0, 1]\times[0, T'])}\leq
    C_{15}.\label{lembdry-1}
        \end{equation}
\end{lem}
\begin{proof}From $(\ref{1.9})_2$, we get
\begin{eqnarray*}
&&\rho\partial_{x}u
=\frac{r^{1-n}}{2c_{1}+c_{2}}\rho^{\gamma-\theta} -(n-1)ur^{-n}
-\frac{\rho^{-\theta}r^{1-n}}{2c_{1}+c_{2}}\int_{x}^{1}\frac{\partial_{t}u}{r^{n-1}}dy\\
&&+\frac{2c_{1}(n-1)}{2c_{1}+c_{2}}\left(r^{-n}u+
\rho^{-\theta}r^{1-n}\int_{x}^{1}\rho^{\theta}\partial_{x}\left(
\frac{u}{r}\right)dy \right)
+\frac{\rho^{-\theta}r^{1-n}}{2c_{1}+c_{2}}\int_{x}^{1}\frac{f}{r^{n-1}}dy.
\end{eqnarray*}
It's easy to show that the first,  second and  fifth terms of the
right side of above equation have bounded $L^\infty$ norm. Similar
to the proof of (\ref{2})$\sim$(\ref{4}), we get the same result
of the fourth term. Now we turn to estimate the third term. Using
 (\ref{2.20}) and $0<\alpha\theta\leq\frac{1}{2}$, we have
    \begin{equation*}
    \left|\frac{\rho^{-\theta}r^{1-n}}{2c_{1}+c_{2}}\int_{x}^{1}\frac{\partial_{t}u}{r^{n-1}}dx\right|
    \leq C\rho^{-\theta}(\int_0^1|\partial_tu|^2dx)^{\frac{1}{2}}
      (1-x)^{\frac{1}{2}}
             \leq C
          (1-x)^{\frac{1}{2}-\alpha\theta}
            \leq C.
    \end{equation*}
From above all, we get (\ref{lembdry}) immediately. From
$(\ref{1.9})_1$, (\ref{rhobdry}) and (\ref{lembdry}), we can
obtain (\ref{lembdry-1}) immediately.
\end{proof}

\section{The proof of Theorem \ref{thm}}\label{sec3}
To construct a weak solution to the initial boundary value problem
(\ref{1.9})$\sim$(\ref{1.12}), we apply the space-discrete
difference scheme method as in \cite{Chen}, which can be described
as follows.

Let $h$ be an increment in $x$,  such that $Nh=1$ for some
$N\in\mathbb{Z}^{+}$, and  $x_{j}=jh$ for $j\in\{0, \cdots, N\}$.
For each integer $N$,  we construct the following time-dependent
functions:
$$(\rho_{j}(t), u_{j}(t), r_{j}(t)), j=0, \cdots, N$$
that form a discrete approximation to $(\rho, u, r)(x_{j}, t)$ for
$j=0, \cdots, N$.

 First,  $(\rho_i(t), u_j(t), r_{j}(t))$,  $i=0, \cdots, N$,
$j=1, \cdots, N$,
        are determined by the following system of $3N+2$
        differential equations:
            \begin{eqnarray}
           &   \dot{\rho}_i=-\rho^2_i\delta(r^{n-1}_iu_i), \label{De1}
          \\&
          \dot{u}_j=r^{n-1}_j\delta
          \sigma_j-2c_{1}(n-1)r^{n-2}_ju_j\delta(\rho^{\theta}_{j-1})+f_j,
          \label{De2}
        \\&
            \displaystyle  r^n_{i+1}=a^n+n\sum^{i}_{k=0}\frac{h}{\rho_k},
               \label{De3}
            \end{eqnarray}
with boundary conditions
    \begin{equation*}
      u_0(t)=0,\  r_0(t)=a,
\end{equation*}
\begin{equation}\label{3.12}
\rho^\gamma_{N}-\rho_{N}^{1+\theta}(2c_1+c_2)\delta(r_{N}^{n-1}u_{N})
+2c_1(n-1)\frac{u_{N+1}}{r_{N+1}}\rho^\theta_{N}=0,
\end{equation}
and initial data
    \begin{equation}
      \rho_j(0)=\frac{1}{h}\int^{jh}_{(j-1)h}\rho_0(y)dy,
    \  u_j(0)=
      \frac{1}{h}\int^{jh}_{(j-1)h}u_0(y)dy,
    \end{equation}
where  $j=1, \cdots, N$,
        \begin{eqnarray}
      &    \rho_0(0)=\rho_{1}(0),
        u_0(0)=0, \ r_0(0)=a,
       \\&
          \displaystyle  r^n_i(0)=a^n+n\sum^{i-1}_{j=0}\frac{h}{\rho_j(0)},
          \ i=1, \cdots, N+1.
            \end{eqnarray}
Here,  $\delta$ is the operator defined by $\delta
w_j=(w_{j+1}-w_j)/h$,  and
    \begin{equation}
      \sigma_j(t)=(2c_{1}+c_{2})\rho^{\theta+1}_{j-1}
      \delta(r^{n-1}_{j-1}u_{j-1})-\rho^{\gamma}_{j-1},
   \               f_j(t)=f(r_j, t).
            \end{equation}

The basic theory of ordinary differential equations guarantees the
local existence and uniqueness of smooth solutions $(\rho_i, u_i,
r_i)$,  $i=0, \cdots, N$, to the system (\ref{De1})-(\ref{3.12})
on some interval $[0, T^h]$,  such that
    $$
    0<\rho_i(t)<\infty,  |u_j(t)|, |r_j(t)|<\infty,\  t\in[0,T^h],
    $$
where $i=0, \cdots,    N$ and $j=0, \cdots,    N+1$. Here, we use
the fact $\displaystyle{\min_{i=0,\ldots,N}}\rho_i(0)>0$.

So, we can let $T^h_{max}$ be the maximal time such that the
smooth solutions exist on $[0, T^h_{max})$ and satisfy
    $$
    0< \rho_i(t)< \infty,  |u_j(t)|, |r_j(t)|<\infty,\
    t\in[0,T_{max}^h),
    $$
where $i=0, \cdots,    N$ and $j=0, \cdots,    N+1$. Our first
goal is to show that $T_{max}^h> T_1$, and the solutions satisfy
    \begin{equation}
    \frac{1}{3}\rho_i(0)\leq \rho_i(t)\leq3\rho_i(0),
   \ i=0, \cdots,    N,\label{local-E3.9}
    \end{equation}
for all $t\in[0,T_1]$, where $T_1>0$ is given in
(\ref{local-E3.8}) and independent of $h$.

Based on the work of Proposition \ref{energy} and Lemma \ref{r},
using similar arguments as that in \cite{Chen, hoff92}, we have
following two lemmas.
\begin{lem}\label{3.1}
Let $(\rho_j(t), u_j(t), r_j(t))$, $j=0, \cdots, N$ be the
solutions to $(\ref{De1})\sim(\ref{3.12})$, then there exists a
positive constant $C(A, \|u_{0}\|_{L^{2}_{x}}, \widetilde{f})$
such that, for all $t\in[0,T^h_{max})$,
\begin{eqnarray}
&&\sum_{j=0}^{N}\left(\frac{1}{2}u_{j}^{2}+\frac{1}{\gamma-1}\rho^{\gamma-1}_j\right)h
       +\int_{0}^{t}\sum_{j=0}^{N}
       \left(
         \left(2c_{1}+c_{2}\right)\rho_j^{\theta+1}\left(\delta\left(r_j^{n-1}u_j\right)\right)^{2}\right.\nonumber\\
          &&+c_1\frac{2(n-1)}{n}\left.\rho_j^{\theta+1}
             (r_j^{n-1}\delta u_j-\frac{u_j}{r_j}\rho_j)^{2}
       \right)h\leq Ce^t.
\end{eqnarray}
\end{lem}
\begin{lem}
The solutions $(\rho_i(t), u_i(t), r_i(t))$, $i=0,\cdots,N$,
satisfy the following identities:
  \begin{eqnarray*}
      &&\qquad\qquad  \partial_{t}r_i(t)=u_i(t), \delta r_i^n(t)=\frac{n}{\rho_i},
\\        &&r_j^{\beta}\delta\rho_j^{\theta}(t)=r_j^{\beta}(0)\delta\rho_j^{\theta}(0)
        -\frac{\theta}{2c_{1}+c_{2}}\left(r_j^{1+\beta-n}(t)u_j(t)-r_j^{1+\beta-n}(0)u_j(0)\right)
     \nonumber\\
        &&\qquad+\frac{\theta}{2c_{1}+c_{2}}
       \left(\int_{0}^{t}
            \left(
            -r_j^{\beta}\delta\rho_j^{\gamma}
            +r_j^{1+\beta+n}f_j+(1+\beta-n)u_j^{2}r_j^{\beta-n})(s)ds
            \right)
       \right),
    \end{eqnarray*}
for all $i=0, \cdots, N$, $j=0, \cdots, N-1$ and
$t\in[0,T^h_{max})$.
\end{lem}

Based on the work of Proposition \ref{rho} and \ref{local-P3.4},
using similar arguments  as that in \cite{Chen, hoff92}, we have
following lemma.

\begin{lem}\label{local-L2.3}
Under the assumptions in Theorem \ref{thm}, we have that for all
$h\in(0, h_0]$($h_0$ is a sufficiently small positive constant ),
there is a $T_{1}>0$ such that, if
    $$
    \frac{1}{3}\rho_i(0)\leq \rho_i(t)\leq 3\rho_i(0),
    $$
for all $i=0,\ldots,N$ and $t\in[0,T']$ where
$T'\in(0,T_{max}^h)\cap(0,T_1]$, then we have
    \begin{eqnarray}
&   0< \dfrac{1}{2}\rho_i(0)\leq\rho_i(t)\leq
2\rho_i(0),\label{local-E2.12}
    \\&
                a\leq r_j(x,t)\leq C,
       \\&
            |u_j(x,t)|\leq C,\label{local-E2.14}
            \end{eqnarray}
for all $i=0,\cdots,N$, $j=0,\cdots,N+1$ and $t\in[0,T']$, where
$C$ is independent of $T'$ and h.
\end{lem}
\begin{rem}
For simplicity of presentation, in Section \ref{sec2}, we
establish some \textit{a priori} estimates in the continuous
version to the initial boundary value problem
(\ref{1.9})$\sim$(\ref{1.12}), so we need $h_0$ is a sufficiently
small positive constant.
\end{rem}

Then, applying the continuation method, we obtain the following
lemma.
\begin{lem}\label{local-L2.4}
Under the assumptions in Lemma \ref{local-L2.3}, we have
(\ref{local-E3.9}) holds for all $t\in [0,T_{max}^h)\cap[0,T_1]$.
\end{lem}
\begin{proof}
  Let $\mathcal{A}=\{T'\in
  [0,T_{max}^h)\cap[0,T_1]\left|\ (\ref{local-E3.9})\ \textrm{ holds for all
  }\ t\in [0,T']
  \right.\}$.

  Since $\min\rho_i(0)>0$, we have $\frac{\rho_i(t)}{\rho_i(0)}\in
  C([0,T_{max}^h))$, $i=0,\ldots,N$, and there exists $T_1^h\in(0,T_{max}^h)$ such that (\ref{local-E3.9}) holds for all $t\in
  [0,T_{1}^h]$. Thus, $\mathcal{A}$ is not empty and relatively
  closed in $[0,T_{max}^h)\cap[0,T_1]$.
    To show that $\mathcal{A}$ is also relatively open in $[0,T_{max}^h)\cap[0,T_1]$, and hence the entire interval,
    it therefore suffices to show that the weaker bound
            $$
            \frac{1}{3}\rho_i(0)\leq \rho_i(t)\leq 3\rho_i(0),
            \ i=0,\ldots,N,\ t\in[0,T']\subset
            [0,T_{max}^h)\cap[0,T_1],
            $$
  implies (\ref{local-E2.12}) holds for all
 $t\in [0,T']$. Thus, from Lemma \ref{local-L2.3}, we have
 $\mathcal{A}=[0,T_{max}^h)\cap[0,T_1]$.
\end{proof}

Then,  applying the reduction to absurdity, we have the following
lemma.
\begin{lem}\label{local-L3.4}
Under the assumptions in Lemma \ref{local-L2.3}, we have
$T_{max}^h> T_1$.
\end{lem}
\begin{proof}
If $T_{max}^h\leq T_1$, from Lemmas
\ref{local-L2.3}$\sim$\ref{local-L2.4}, we have
(\ref{local-E2.12})$\sim$(\ref{local-E2.14}) hold on
$t\in[0,T_{max}^h)$. Thus, we can extend the existence interval
$[0,T_{max}^h)$ to $[0,T_{max}^h]$.  It contradicts the definition
of $T_{max}^h$.
\end{proof}

Based on the work of Lemmas \ref{u4m}$\sim$\ref{lem}, using
similar arguments as that in \cite{Chen, hoff92}, we have
following lemma.

\begin{lem}\label{local-L3.5}
Under the assumptions  in Lemma \ref{local-L2.3}, we have that
$(\rho_j(t), u_j(t), r_j(t))$, $j=0,\cdots,N$,  satisfy
    \begin{eqnarray}\label{3.17}
    \sum_{j=0}^{N}u_j^{4m}(t)h+\int_0^t\sum_{j=0}^{N}\left[\rho_j^{\theta+1}(s)u_j^{4m-2}(s)
    \delta(r_j^{n-1}(s)u_j(s))^{2}\right.\nonumber\\
    \left.+\rho_j^{\theta-1}(s)r_j^{-2}(s)u_j^{4m}(s)\right]hds\leq C,
    \end{eqnarray}
        \begin{equation}
        \sum_{j=0}^{N-1}(1-jh)^{\alpha_{0}}(\delta\rho_j^{\theta}(t))^{2}h+
    \sum_{j=0}^{N}(\partial_{t}u_j)^{2}h+
            \sum_{j=0}^{N}|\delta u_j(t)|^{\lambda_0}h\leq C,
            \end{equation}
    \begin{equation}
    \sum_{j=0}^{N}|\frac{d}{dt}\rho_j(t)|^{2}h+
                \int^{T_1}_0\sum_{j=0}^{N}
                \left|\frac{d}{dt}\left(\rho_{j}^{\theta+1}\delta(r_{j}^{n-1}u_{j})\right)(t)
                \right|^{2}hdt\leq        C,
            \end{equation}
    \begin{equation}
    \|\rho_{i}^{\theta+1}\delta(r_{i}^{n-1}u_{i})(t)\|_{L^{\infty}([0,
    T_{1}])}+\|u_i(t)\|_{L^{\infty}([0, T_{1}])}\leq C,
    \end{equation}
        \begin{equation}
        \sum_{j=0}^{N-1}|\delta(\rho_{j}^{\theta+1}\delta(r_{j}^{n-1}u_{j}))(t)|h+
        \sum_{j=0}^{N-1}|\delta\rho_{j}(t)|h\leq C,
        \end{equation}
            \begin{equation}
              |r_i(t)|+\sum_{j=0}^{N}|\delta
              r_j|^{\lambda_0}h+|\frac{d}{dt}r_i|\leq C,
            \end{equation}
  for  $i=0, \cdots, N$,  $0\leq s\leq t\leq
T_{1}$. Furthermore, if $\alpha\theta\leq\frac{1}{2}$, we have
        \begin{equation}
            \|\rho_i\delta (r^{n-1}_iu_i)\|_{L^\infty([0,T_1])}+
            \|\rho_i^{-1}(0)\frac{d\rho_i(t)}{dt}\|_{L^\infty([0,T_1])}\leq
            C,\ i=0, \cdots, N.\label{local-E3.26}
        \end{equation}
\end{lem}
Now, we turn to prove \textbf{Theorem \ref{thm}}.

With $(\rho_i, u_i, r_i)$, $i=0,\cdots,N$, we can define our
approximate solutions $(\rho^N$,$ u^N$, $r^N)(x, t)$ for the
system (\ref{1.9})$\sim$(\ref{1.12}). For each fixed $N$ and
$t\in[0, T_1]$,  we define piecewise linear continuous functions
$(\rho^N, u^N, r^N)(x, t)$ with respect to $x$ as follows: when
$x\in([xN]h, ([xN]+1)h]$,
\begin{equation}
\rho^N(x,
t)=\rho_{[xN]}(t)+(xN-[xN])(\rho_{[xN]+1}(t)-\rho_{[xN]}(t)),
\label{E3.12}
\end{equation}
    \begin{equation}\label{3.14}
    u^N(x, t)=u_{[xN]}(t)+(xN-[xN])(u_{[xN]+1}(t)-u_{[xN]}(t)),
    \end{equation}
\begin{equation}\label{3.15}
r^N(x,
t)=\left(r^n_{[xN]}(t)+(xN-[xN])(r^n_{[xN]+1}(t)-r^n_{[xN]}(t))
\right)^{1/n}.
\end{equation}
We have for $jh\leq x\leq(j+1)h$
\begin{eqnarray}
\partial_{x}u^{N}(x, t)=\frac{u_{j+1}(t)-u_{j}(t)}{h}.
\end{eqnarray}
We also introduce the corresponding step functions:
\begin{equation}
(\rho_h, u_h, r_h)(x, t)=(\rho_{[xN]}, u_{[xN]},
r_{[xN]})(t),\indent x\in([xN]h, ([xN]+1)h].
\end{equation}

Using Lemmas \ref{3.1}-\ref{local-L3.5}, we have, for all $h\in(0,
h_0]$
    $$
  \sup_{t\in[0,T_1]}\left(
  \|u^N(\cdot,t)\|_{L^\infty\cap W^{1,\lambda_0}}+
  \|u^N_t(\cdot,t)\|_{L^2}+\|\partial_t(\rho_h)(\cdot,t)\|_{L^2}
  \right)\leq C,
    $$
        $$
        \frac{1}{3}\rho^N(x,0)\leq \rho^N(x,t)\leq 3\rho^N(x,0),
       \ \forall\ (x,t)\in[0,1]\times[0,T_1],
        $$
        $$
      \sup_{t\in[0,T_1]}\left(
      \int^1_0(1-x)^{\alpha_0}((\rho^N)^\theta)_x^2dx+TV_{[0,1]}(\rho^N)
      \right)\leq C,
        $$
    $$
    \sup_{t\in[0,T_1]}\left(
        \|(\rho_h^{1+\theta}\partial_x((r^N)^{n-1}u^N))\|_{L^\infty_x}+
        TV_{[0,1]}(\rho_h^{1+\theta}\partial_x((r^N)^{n-1}u^N))
        \right)\leq C,
    $$
        $$
        \int^{T_1}_0\|\partial_t(\rho_h^{1+\theta}\partial_x((r^N)^{n-1}u^N))(\cdot,t)\|^2_{L^2}dt\leq
        C,
        $$
    $$
    \sup_{t\in[0,T_1]}\left(
    \|r^N(\cdot,t)\|_{L^\infty\cap
    W^{1,\lambda_0}}+\|\partial_tr^N(\cdot,t)\|_{L^\infty}
    \right)\leq C,
    $$
where $TV(g)$ is the total variation of $g$, $x\in[0,1]$ and
$t\in[0,T_1]$.

With the above estimates, using Helly's theorem and similar
arguments as that in Section 9 of \cite{Chen}, we can get the
following compactness of the approximate solutions (If necessary,
we can choose the subsequence.):
\begin{eqnarray}\label{compactness}
&&(\rho_h,\rho^N, u^N,
r^N)(x, t)\rightarrow(\rho,\rho, u, r)(x, t),\indent a.e.\nonumber\\
&& (\rho_h)^{1+\theta}\partial_x((r^N)^{n-1}u^N)(x, t)\rightarrow
\rho^{1+\theta}\partial_x(r^{n-1}u)(x, t),\indent a.e.\nonumber\\
&&(\partial_x r^N, \partial_x u^N)(x,t)\rightharpoonup (\partial_x
r,\partial_x u)(x,t),\indent \textrm{weakly in}\
L^{\lambda_0}([0,1]\times[0,T_1]), \nonumber\\
&&\partial_t\rho^N(x, t)\rightharpoonup\partial_t\rho(x,
t),\indent \textrm{weakly in}\ L^2([0,1]\times[0,T_1]),
\end{eqnarray}
when $N\rightarrow\infty$, where $r(x,
t)=(a^n+n\int_0^x\frac{dy}{\rho(y, t)})^{\frac{1}{n}}$.

For any given test function $\psi(x, t)\in C_0^\infty((0,
1)\times[0, T_1])$, we choose $N=\frac{1}{h}$ is large enough such
that the support of the test function $\psi$ is away  enough from
the boundaries, that is, supp $\psi\subset((h, 1-h)\times[0,
T_1])=((\frac{1}{N}, 1-\frac{1}{N})\times[0, T_1])$.

Define
\begin{equation}
\psi_i( t)=\psi_h(x, t)=\psi([xN]h, t),\indent ih\leq x< (i+1)h.
\end{equation}
 We can
see that $\psi_i(t)=0$ for $i=0, N-1, N$.

Multiplying (\ref{De1}) by $\psi_i$, summing it up for $i=0,
\cdots, N$, and integrating it over $[0, T_1]$, we get
\begin{eqnarray}\label{pass1}
0&=&\int_0^{T_1}\sum_{i=0}^{N}
   \left(
   \psi_i\partial_t\rho_i+\psi_i\rho_i^2\delta(r_i^{n-1}u_i)
   \right)hdt\nonumber\\
&=&O(h)+\int_0^{T_1}\int_0^1\left(
         \psi_h\partial_t\rho_h+\psi_h(\rho_h)^2\partial_x((r^N)^{n-1}u^N)
                            \right)dxdt\nonumber\\
&=&O(h)+\int_0^{T_1}\int_0^1\left(
         \psi\partial_t\rho_h+\psi(\rho_h)^2\partial_x((r^N)^{n-1}u^N)
                            \right)dxdt.
\end{eqnarray}
Using (\ref{compactness}), we can pass the limit in (\ref{pass1})
to obtain the first equation in (\ref{1.9}).

For any given test function $\phi(x, t)\in C_0^\infty((0,
1]\times[0, T_1))$, we choose $N=\frac{1}{h}$ is large enough such
that the support of the test function $\phi$ is away  enough from
the fixed boundary, that is, $\phi|_{x\in[0,h]}=0$.

Define
\begin{equation}
\phi_i(  t)=\phi_h(x, t)=\phi([xN]h, t),
\end{equation}
    \begin{equation}
    \phi^N(x,t)=\phi_{[xN]}(t)+(xN-[xN])
    (\phi_{[xN]+1}(t)- \phi_{[xN]}(t)),
    \end{equation}
when $ih\leq x< (i+1)h$.
 We can
see that $\phi_i(t)=0$ for $i=0, N-1, N$.

Multiplying (\ref{De2}) by $\phi_j$, summing it up for $j=1,
\cdots, N$, and integrating it over $[0, T_1]$, we get
\begin{eqnarray}\label{pass2}
0&=&\int_0^{T_1}\sum_{j=1}^{N}\phi_j\partial_t u_j hdt
-\int_0^{T_1}\sum_{j=1}^{N}\phi_jr_j^{n-1}\delta\sigma_jhdt\nonumber\\
&&+2c_1(n-1)\int_0^{T_1}\sum_{j=1}^{N}\phi_jr_j^{n-2}u_j\delta(\rho_{j-1}^{\theta})hdt
-\int_0^{T_1}\sum_{j=1}^{N}\phi_jf_jhdt\nonumber\\
&=&O(h)-\int_0^1\phi_hu^N(x,
0)dx-\int_0^{T_1}\int_0^1\partial_t\phi_hu^N(x, t)dxdt\nonumber\\
&&+\int_0^{T_1}\int_0^1(2c_1+c_2)\partial_x(\phi^Nr^N)\rho_h^{\theta+1}\partial_x((r^N)^{n-1}u^N)dxdt\nonumber\\
&&-\int_0^{T_1}\int_0^1\partial_x(\phi^Nr^N)\rho_h^{\gamma}dxdt
-\int_0^{T_1}\int_0^12c_1(n-1)
 \left(
 \partial_x\phi^N(r^{N})^{n-2}u^N\right.\nonumber\\
 &&\left.+(n-2)\phi_h(r^N)^{n-3}\partial_xr^Nu^N
 +\phi_h(r^N)^{n-2}\partial_xu^N
 \right)\rho_h^{\theta}dxdt\nonumber\\
&& -\int_0^{T_1}\int_0^1\phi_hf_hdxdt,
\end{eqnarray}
where $f_h(x, t)=f_{[xN]}(t)$, $x\in([xN]h, ([xN]+1)h]$. Using
(\ref{compactness}), we can pass the limit in (\ref{pass2}) to
obtain (\ref{1.9_2}). Thus we complete the proof of Theorem
\ref{thm}.

\begin{rem}
When $\alpha\theta\leq\frac{1}{2}$,  from (\ref{local-E3.26}), we
have
    $$
    \|\rho^N((r^N)^{n-1}u^N)_x\|_{L^\infty([0,1]\times[0,T_1])}+\|(\rho^N(x,0))^{-1}(\rho^N)_t
    \|_{L^\infty([0,1]\times[0,T_1])}\leq C.
    $$
Thus, the limit function $(\rho,u,r)$ satisfies
    $$
    \|\rho(r^{n-1}u)_x\|_{L^\infty([0,1]\times[0,T_1])}+\|\rho_0^{-1}\rho_t
    \|_{L^\infty([0,1]\times[0,T_1])}\leq C,
    $$
and
    $$
    (1-x)^{-\alpha}\rho\in C([0,T_1];L^\infty([0,1])).
    $$
\end{rem}

\begin{rem}
There is another method to prove  Theorem \ref{thm}. At first, we
consider system
    \begin{equation}\label{local-E1.9}
\left \{ \begin{array}{l}
\partial_{t}\rho=-\rho^{2}\partial_{x}(r^{n-1}u) \\
\partial_{t}u=r^{n-1}\partial_{x}\left(\rho(\lambda+2\mu)\partial_{x}(r^{n-1}u)-P\right)
              -2(n-1)r^{n-2}u\partial_{x}\mu+f(r(x, t), t)\\
r^{n}(x, t)=a^{n}+n\int_{0}^{x}\rho^{-1}(y, t)dy
\end{array} \right.
\end{equation}
with the initial data
    \begin{equation}
    (\rho, u)|_{t=0}=(\rho^\varepsilon_{0}, u^\varepsilon_{0})(x),
    r|_{t=0}=r^\varepsilon_{0}(x)
    =\left(a^{n}+n\int_{0}^{x}(\rho^\varepsilon_{0})^{-1}(y)dy\right)^{\frac{1}{n}},
    \end{equation}
where $\rho_0^\varepsilon(>\varepsilon), u^\varepsilon_0,
r_0^\varepsilon$ converge to $\rho_0, u_0, r_0$ in some suitable
spaces as $\varepsilon$ goes to 0, and the boundary conditions:
    \begin{equation}
    u(0, t)=0,
    \end{equation}
        \begin{equation}
    \left.\left\{\rho(\lambda+2\mu)\partial_{x}(r^{n-1}u)-P
              -2(n-1)r^{-1}u\mu\right\}\right|_{x=1}=0.\label{local-E1.12}
    \end{equation}
Using similar arguments as that in Sections
\ref{sec2}$\sim$\ref{sec3}, we can obtain the existence of the
weak solution $(\rho^\varepsilon,u^\varepsilon,r^\varepsilon)$ to
the system (\ref{local-E1.9})$\sim$(\ref{local-E1.12}), and some
uniform estimates of the solution. Letting
$\varepsilon\rightarrow0$, we can prove that  the limit function
$(\rho,u,r)$ is the weak solution to the system
(\ref{1.9})$\sim$(\ref{1.12}).
\end{rem}
\section{Continuous dependence}\label{sec4}
In this section,  we will prove Theorem \ref{Continuous
Dependence}, applying the energy method. Let $(\rho_1, u_1,
r_1)(x, t)$ and $(\rho_2, u_2, r_2)(x, t)$ be two solutions in
Theorem \ref{thm} corresponding to the initial data $(\rho_{01},
u_{01}, r_{01})(x, t)$ and $(\rho_{02}, u_{02}, r_{02})(x, t)$,
respectively. Then we have, $i=1, 2$, $(x, t)\in[0, 1]\times[0,
T_{1}]$,
    \begin{equation}
      \frac{A}{3}(1-x)^{\alpha}\leq \rho_i(x, t)\leq 3B(1-x)^{\alpha},
      \ |u_i(x, t)|\leq C,
      \ a\leq r_i\leq C.
      \label{4.1}
    \end{equation}
From Lemma \ref{lem}, we can easily get
    \begin{equation}
    \|\rho_i\partial_x(r_i^{n-1}u_i)\|_{L^{\infty}([0, 1]\times[0,
    T_{1}])}\leq C.
     \label{4.2}
    \end{equation}

For simplicity,  we may assume that $(\rho_1, u_1,r_1)(x, t)$ and
$(\rho_2, u_2,r_2)(x, t)$ are suitably smooth since the following
estimates are valid for the solutions with the regularities
indicated in Theorem $\ref{thm}$ by using the Friedrichs
mollifier.

Let
    \begin{eqnarray*}
      \varrho=\rho_1-\rho_2,
      \ w=u_1-u_2,
      \ \mathcal{R}=r_1-r_2.\\
      \varrho_0=\rho_{01}-\rho_{02},
      \ w_0=u_{01}-u_{02},
      \ \mathcal{R}_0=r_{01}-r_{02}.
    \end{eqnarray*}
From $(\ref{1.9})_1$, (\ref{4.1})$\sim$(\ref{4.2}), Lemma \ref{r}
and Young's inequality, we have
    \begin{eqnarray}
      &&\frac{d}{dt}\int^1_0\rho_1^\theta\rho_2^{-1}\mathcal{R}^2dx\nonumber\\
            &=&\int^1_02\rho_1^\theta\rho_2^{-1}\mathcal{R}\mathcal{R}_t
             -\theta\rho_1^{\theta+1}\partial_x(r_1^{n-1}u_1)\rho_2^{-1}\mathcal{R}^{2}
             +\rho_1^{\theta}\partial_x(r_2^{n-1}u_2)\mathcal{R}^{2}dx\nonumber\\
        &\leq&C_{\varepsilon}\int^1_0\rho_1^\theta\rho_2^{-1}\mathcal{R}^2dx
             +\varepsilon\int_0^1\frac{w^2}{r_1^2}\rho_1^{\theta-1}dx,
    \end{eqnarray}
where $\varepsilon>0$ is chosen later.

From (\ref{1.9})$_1$, we have
\begin{eqnarray}
  \partial_t \varrho&=&-\rho^2_1r^{n-1}_1\partial_xw
        -\varrho(\partial_xu_2r^{n-1}_1(\rho_1+\rho_2)+(n-1)\frac{u_2}{r_2})
        \nonumber\\
  &&-(n-1)\frac{\rho_1}{r_1}w-\rho_2^2\partial_xu_2(r^{n-1}_1-r^{n-1}_2)
  +(n-1)\frac{u_2\rho_1}{r_1r_2}\mathcal{R},
\end{eqnarray}
and
    \begin{eqnarray}
    &&\frac{d}{dt}\int^1_0\rho_1^{1-\theta}\rho_2^{2\theta-4}\varrho^2dx\nonumber\\
    &=&\int^1_02\rho_1^{1-\theta}\rho_2^{2\theta-4}\varrho\partial_t\varrho
       -(1-\theta)\rho_1^{2-\theta}\partial_x(r_1^{n-1}u_1)\rho_2^{2\theta-4}\varrho^2\nonumber\\
    &&-(2\theta-4)\rho_1^{1-\theta}\partial_x(r_2^{n-1}u_2)\rho_2^{2\theta-3}\varrho^2dx\nonumber\\
    &\leq&2\int^1_0\rho_1^{1-\theta}\rho_2^{2\theta-4}\varrho\partial_t\varrho dx
     +C\int_0^1\rho_1^{1-\theta}\rho_2^{2\theta-4}\varrho^2dx\nonumber\\
    &=&2\int_0^1-\rho^{3-\theta}_1\rho_2^{2\theta-4}\varrho r^{n-1}_1\partial_xw-
        \rho_1^{1-\theta}\rho_2^{2\theta-4}\varrho^2
        (\partial_xu_2r^{n-1}_1(\rho_1+\rho_2)+(n-1)\frac{u_2}{r_2})
        \nonumber\\
  &&-(n-1)\frac{\rho_1^{2-\theta}\rho_2^{2\theta-4}\varrho}{ r_1}w
  -\rho_1^{1-\theta}\rho_2^{2\theta-2}\varrho\partial_xu_2(r^{n-1}_1-r^{n-1}_2)\nonumber\\
  &&+(n-1)\frac{u_2\rho_1^{2-\theta}\rho_2^{2\theta-4}\varrho}{r_1r_2}\mathcal{R}dx
  +C\int_0^1\rho_1^{1-\theta}\rho_2^{2\theta-4}\varrho^2dx\nonumber\\
    &\leq&C_{\varepsilon}\int_0^1\rho_1^{1-\theta}\rho_2^{2\theta-4}\varrho^2 dx
         +C_{\varepsilon}\int_0^1\rho_1^{-1}\rho_2^{\theta}\mathcal{R}^2dx
         +\varepsilon\int_0^1\frac{w^2}{r_1^2}\rho_1^{\theta-1}dx\nonumber\\
         &&+\varepsilon\int_0^1\rho_1^{\theta+1}r_1^{2n-2}(\partial_xw)^2dx.
    \end{eqnarray}
\begin{rem}
We can use the weighted function
$\rho_1^{l_1}\rho_2^{l_2}(\rho_1-\rho_2)^2$, where $l_1,l_2$
satisfy $l_1+l_2=\theta-3$. For simple, we choose $l_1=1-\theta$,
$l_2=2\theta-4$.
\end{rem}
We only give a part of the proof of the last inequality. The rests
are the same.
\begin{eqnarray*}
&&\int_0^1-\rho^{3-\theta}_1\rho_2^{2\theta-4}\varrho
r^{n-1}_1\partial_xwdx\leq
R^{n-1}\int_0^1\rho_1^{\frac{\theta+1}{2}}\partial_xw
               \centerdot\rho_1^{\frac{1-\theta}{2}}\rho_2^{\theta-2}\varrho
               \centerdot\rho_1^{2-\theta}\rho_2^{\theta-2}dx\nonumber\\
&\leq&R^{n-1}(\frac{9B}{A})^{2-\theta}\int_0^1|\rho_1^{\frac{\theta+1}{2}}\partial_xw
               \centerdot\rho_1^{\frac{1-\theta}{2}}\rho_2^{\theta-2}\varrho|dx\nonumber\\
&\leq&\varepsilon\int_0^1\rho_1^{\theta+1}r_1^{2n-2}(\partial_xw)^2dx
      +C_{\varepsilon}\!\!\int_0^1\rho_1^{1-\theta}\rho_2^{2\theta-4}\varrho^2dx.
\end{eqnarray*}
From equation (\ref{1.9})$_2$ and boundary conditions
(\ref{1.10})-(\ref{1.12}),  we get
    \begin{eqnarray}
     \frac{d}{dt}\int^1_0\frac{1}{2}w^2(x, t)dx
            &=&\int^1_0\left\{-(\lambda(\rho_1)+2\mu(\rho_1))\rho_1\partial_x(r^{n-1}_1u_1)
            \partial_x(r^{n-1}_1 w)\right.
            \nonumber\\
            &&\left.+(\lambda(\rho_2)+2\mu(\rho_2))\rho_2\partial_x(r^{n-1}_2u_2)
            \partial_x(r^{n-1}_2 w)
            \right\}dx\nonumber\\
                &&+\int^1_0\left\{P(\rho_1)\partial_x(r^{n-1}_1w)-P(\rho_2)\partial(r^{n-1}_2w)
                \right\}dx\nonumber\\
      &&+2(n\!-\!1)\int^1_0\left\{\mu(\rho_1)\partial_x(r^{n-2}_1u_1w)
      \!-\!\mu(\rho_2)\partial_x(r^{n-2}_2u_2w)
      \right\}dx\nonumber\\
        &&+\int^1_0(f(r_1, t)-f(r_2, t))wdx:=U_1+U_2+U_3+U_4.
    \end{eqnarray}
Using the similar argument as that in Proposition \ref{energy},
Cauchy-Schwartz inequality and (\ref{4.1})$\sim$(\ref{4.2}), we
have
    \begin{eqnarray}
      U_1+U_3
      &\leq&-\frac{1}{2}\int^1_0\{(c_2+\frac{2}{n}c_1)\rho_1^{1+\theta}[\partial_x(r_1^{n-1}w)]^2\nonumber\\
      &&+\frac{2c_1(n-1)}{n}
      \rho_1^{\theta+1}[r_1^{n-1}\partial_x w-\frac{w}{r_1\rho_1}]^2\}dx\nonumber\\
      &&+C\int_0^1(\rho_1^{1-\theta}\rho_2^{2\theta-4}\varrho^2
                    +\rho_1^{\theta}\rho_2^{-1}\mathcal{R}^2)dx,
   \\
          U_2 &\leq& \varepsilon\int^1_0\!\rho_1^{1+\theta}r^{2n-2}_1
    [\partial_xw]^2dx + \varepsilon\int_0^1\!\rho_1^{\theta-1}w^2dx\nonumber\\
               &&  + C_{\varepsilon}\int^1_0\! (\rho_1^{1-\theta}\rho_2^{2\theta-4}\varrho^2
                + \rho_1^{\theta}\rho_2^{-1}\mathcal{R}^2)dx,
        \end{eqnarray}
and
    \begin{eqnarray}
      U_4&=&\int^1_0(f(r_1, t)-f(r_2, t))wdx\\
      &\leq&
      C\int_0^1|\mathcal{R}w|dx
      \leq
      C_{\varepsilon}\int_0^1\rho_1^{\theta-1}\mathcal{R}^2dx
        +\varepsilon\int_0^1\rho_1^{1-\theta}w^2dx\nonumber\\
      &\leq&
      C_{\varepsilon}\int_0^1\rho_1^{\theta-1}\mathcal{R}^2dx
        +\varepsilon C\|w\|_{L^{\infty}}^2\nonumber\\
      &\leq&
      C_{\varepsilon}\int_0^1\rho_1^{\theta-1}\mathcal{R}^2dx
        +\varepsilon C\int_0^1\rho_1^{\theta+1}r_1^{2n-2}(\partial_xw)^2dx
        \int_0^1\rho_1^{-\theta-1}r_1^{-2n+2}dx\nonumber\\
      &\leq&
      C_{\varepsilon}\int_0^1\rho_1^{\theta-1}\mathcal{R}^2dx
        +\varepsilon C\int_0^1\rho_1^{\theta+1}r_1^{2n-2}(\partial_xw)^2dx.    \end{eqnarray}
Here,  we use the fact that $f(r, t)\in C^1([a, +\infty)\times[0,
+\infty))$.
 Choosing a sufficiently small positive constant $\varepsilon$, we obtain
    \begin{eqnarray*}
      &&\!\!\!\!\!\frac{d}{dt}\int^1_0\!\!
        \left(
         w^2
        \! +\!\rho_1^{1-\theta}\rho_2^{2\theta-4}\varrho^2
        \! +\!\rho_1^{\theta}\rho_2^{-1}\mathcal{R}^2
        \right)dx\!+\!C_{16}\int^1_0\!\!\rho_1^\theta\left(\rho_1r^{2n-2}_1\left(\partial_xw\right)^2
        \!+\!      \frac{w^2}{r^2_1\rho_1}\right)dx\\
        &&\leq C_{17}\int^1_0
        \left(
        w^2
        +\rho_1^{1-\theta}\rho_2^{2\theta-4}\varrho^2
        +\rho_1^{\theta}\rho_2^{-1}\mathcal{R}^2
        \right)dx,
    \end{eqnarray*}
where $C_{16}, C_{17}$ are two positive constants dependent on
$\varepsilon$.

Using Gronwall's inequality,  we have for any $t\in[0, T_{1}]$,
\begin{eqnarray*}
 \int^1_0
        \left(
         w^2
         \!+\!\rho_{1}^{1-\theta}\rho_{2}^{2\theta-4}\varrho^2
         \!+\!\rho_{1}^{\theta}\rho_{2}^{-1}\mathcal{R}^2
        \right)dx
 \leq  Ce^{Ct}\int^1_0
        \left(
         w_0^2
         \!+\!\rho_{01}^{1-\theta}\rho_{02}^{2\theta-4}\varrho_0^2
        \! +\!\rho_{01}^{\theta}\rho_{02}^{-1}\mathcal{R}_0^2
        \right)dx.
\end{eqnarray*}
 Then,
we finish the proof of Theorem \ref{Continuous Dependence}.

\section{Appendix}\label{sec5}
\subsection{Proof of Lemma \ref{u4m}}
\begin{proof} We apply the inductive method to prove this lemma.

First we consider the case of $k=1$. From Proposition
\ref{energy}, we obtain (\ref{u4m1}) with $k=1$ immediately.

Assume (\ref{u4m1}) holds for $k=l-1$, i.e.,
\begin{equation*}
\int_0^1u^{2(l-1)}dx+\int_0^t\int_0^1\rho^{\theta+1}u^{2(l-1)-2}r^{2n-2}(\partial_xu)^2dxds\leq
C.
\end{equation*}
Thus, using $\alpha<\frac{1}{1+\theta}$, we have
\begin{eqnarray}\label{*}
&&\int_0^t\|u^{l-1}\|_{L^{\infty}}^2ds
\leq\int_0^t\left(\int_0^1|\partial_x(u^{l-1})|dx\right)^2 ds\nonumber\\
&\leq&\int_0^t\left(\int_0^1\rho^{\theta+1}\left(\partial_x(u^{l-1})\right)^2dx\right)\left(\int_0^1\rho^{-\theta-1}dx\right)ds
\leq C.
\end{eqnarray}

Now we need to prove (\ref{u4m1}) holds for $k=l$. Multiplying
$(\ref{1.9})_{2}$ by $u^{2l-1}$ and integrating it over $x$ from 0
to 1,  using the boundary conditions
$(\ref{1.11})\sim(\ref{1.12})$, we have
\begin{eqnarray}\label{2.9}
\frac{d}{dt}\int_{0}^{1}\frac{1}{2l}u^{2l}(x, t)dx
&=&-(2c_{1}+c_{2})\int_{0}^{1}\rho^{\theta+1}\partial_{x}(r^{n-1}u)
   \partial_{x}(r^{n-1}u^{2l-1}) dx\nonumber\\
&&+\int_{0}^{1}\!\rho^{\gamma}\partial_{x}(r^{n\!-\!1}u^{2l-1}) dx
\!+\!2c_{1}(n\!-\!1)\int_{0}^{1}\!\rho^{\theta}\partial_{x}(r^{n-2}u^{2l}) dx\nonumber\\
&&+\int_{0}^{1}fu^{2l-1} dx\nonumber\\
&:=&G_{1}+G_{2}+G_{3}+G_{4}.
\end{eqnarray}
Set
$$A_{1}^{2}=\rho^{\theta+1}u^{2l-2}r^{2n-2}(\partial_xu)^2\geq0, \quad
A_{2}^2=\rho^{\theta-1}r^{-2}u^{2l}\geq0.$$ Thus
\begin{eqnarray*}
G_{1}+G_{3}&=&-(2l-1)(2c_{1}+c_{2})\int_{0}^{1}A_{1}^{2}dx
-2(n-1)lc_2\int_0^1A_{1}A_{2}dx\\
&&+\left(2(n-1)(n-2)c_1-(n-1)^2(2c_1+c_2)\right)\int_0^1A_2^2dx.
\end{eqnarray*}
Then, from (\ref{2.9}) and  Young's inequality, we get
\begin{eqnarray}\label{5.3}
&&\frac{d}{dt}\int_0^1\frac{1}{2l}u^{2l}dx+(2l-1)(2c_{1}+c_{2})\int_{0}^{1}A_{1}^{2}dx\nonumber\\
&\leq&\varepsilon\int_0^1A_1^2dx+C_{\varepsilon}\int_0^1A_2^2dx
      +\int_0^1\rho^{\gamma}\partial_x(r^{n-1}u^{2l-1})dx
      +\int_0^1fu^{2l-1}dx.
\end{eqnarray}
From (\ref{upalpha}), (\ref{rhobdry}), (\ref{*}), Young's
inequality and $r\geq a$, we obtain
\begin{eqnarray*}
\int_0^1A_2^2dx&\leq&\int_0^1C(1-x)^{\alpha(\theta-1)}a^{-2}u^{2l}dx\\
&\leq&C\int_0^1(1-x)^{\alpha(\theta-1)}u^{2l-2}dx+C\int_0^1u^{2l+2}dx\\
&\leq&C\|u^{2l-2}\|_{L^{\infty}}\int_0^1(1-x)^{\alpha(\theta-1)}dx+C\|u^2\|_{L^{\infty}}\int_0^1u^{2l}dx\\
&\leq&C(\|u^{l-1}\|_{L^{\infty}}^2+\|u\|_{L^{\infty}}^2\int_0^1u^{2l}dx),
\end{eqnarray*}
\begin{equation*}
\int_0^1\rho^{\gamma}\partial_x(r^{n-1}u^{2l-1})dx\leq
\varepsilon\int_0^1A_1^2dx+C_{\varepsilon}+C_{\varepsilon}\int_0^1u^{2l}dx,
\end{equation*}
and
\begin{equation*}
\int_0^1fu^{2l-1}dx\leq C+C\int_0^1u^{2l}dx.
\end{equation*}
Inserting above inequalities to (\ref{5.3}) and choosing
$\varepsilon=\frac{(2l-1)(2c_1+c_2)-1}{2}$, we get that there
exist two positive constants $C_{18}$ and $C_{19}$ such that
\begin{eqnarray*}
&&\frac{d}{dt}\int_0^1u^{2l}dx+\int_0^1\rho^{\theta+1}u^{2l-2}r^{2n-2}(\partial_xu)^2dx\\
&\leq&
C_{18}+C_{19}(\|u^{l-1}\|_{L^{\infty}}^2+\|u\|_{L^{\infty}}^2\int_0^1u^{2l}dx).
\end{eqnarray*}
Then, using Gronwall's inequality, (A3) and (\ref{*}), we get
\begin{eqnarray*}
&&\int_{0}^{1}u^{2k}(x,
t)dx+\int_{0}^{t}\int_{0}^{1}\left(\rho^{\theta+1}u^{2k-2}r^{2n-2}(\partial_xu)^2\right)(x,s)dxds\\
&\leq&(C_{18}+C_{19}\int_0^t\|u^{l-1}\|_{L^{\infty}}^2ds)e^{C_{19}\int_0^t\|u\|_{L^{\infty}}^2ds}\leq
C_{2},\indent 0\leq t\leq T'.
\end{eqnarray*}
So (\ref{u4m1}) holds for $k=l$. From all above, we get the result
immediately.
\end{proof}
\subsection{Proof of Lemma \ref{ul2}}
\begin{proof}\ From (\ref{1.9}),  we get
    $$
\partial_{t}(u-u_{0})
=r^{n-1}\partial_{x}\left((2c_{1}+c_{2})\rho^{\theta+1}
\partial_{x}(r^{n-1}u)-\rho^{\gamma}\right)-2c_{1}(n-1)r^{n-2}u\partial_{x}\rho^{\theta}+f.
    $$
Multiply the above equation by $u-u_{0}$ and integrating it over
$x$ from 0 to 1, similar to the proof of proposition \ref{energy},
and using the inequality $ab\leq\frac{1}{p}a^{p}+\frac{1}{q}b^{q}$
, where $\frac{1}{p}+\frac{1}{q}=1, p, q>1, a, b\geq0$,  we have
\begin{eqnarray*}
&&\frac{d}{dt}\int_{0}^{1}\frac{1}{2}(u-u_{0})^{2}dx\\
&&+\int_{0}^{1}(c_{2}+\frac{2}{n}c_1)\rho^{\theta+1}\left(\partial_{x}(r^{n-1}u)\right)^{2}
 +\frac{2(n-1)}{n}c_{1}\rho^{\theta+1}\left(r^{n-1}\partial_{x}u-\frac{u}{r\rho}\right)^{2}dx\\
&=&\int_{0}^{1}(2c_{1}+c_{2})\rho^{\theta+1}\partial_{x}(r^{n-1}u_{0})
    \partial_{x}(r^{n-1}u)dx\\
&&-\int_{0}^{1}2c_{1}(n-1)\left((n-2)r^{-2}\rho^{-1}uu_{0}
  +r^{n-2}\partial_{x}uu_{0}+r^{n-2}u(u_{0})_x\right)\rho^{\theta}dx\\
&&+\int_{0}^{1}\partial_{x}(r^{n-1}u)\rho^{\gamma}dx
  -\int_{0}^{1}\partial_{x}(r^{n-1}u_{0})\rho^{\gamma}dx
+\int_{0}^{1}(u-u_{0})fdx\\
&\leq&\frac{\epsilon}{2}\int_{0}^{1}\rho^{\theta+1}\left(\partial_{x}(r^{n-1}u)\right)^{2}dx
   +\frac{\epsilon}{2}\int_{0}^{1}\rho^{\theta-1}u^{2}r^{-2}dx
   +\frac{\epsilon}{2}\int_{0}^{1}\rho^{\theta+1}(\partial_xu)^{2}r^{2n-2}dx\\
&&
   +C_{\epsilon}\int_{0}^{1}\rho^{\theta+1}\left(\partial_{x}(r^{n-1}u_{0})\right)^{2}dx
   +C_{\epsilon}\int_{0}^{1}\rho^{\theta-1}u_{0}^{2}r^{-2}dx\\
&&
   +C_{\epsilon}\int_{0}^{1}\rho^{\theta-1}r^{-2}u_{0}^{2}dx
   +C_{\epsilon}\int_{0}^{1}\rho^{\theta+1}r^{2n-2}((u_{0})_x)^{2}dx
   +C_{\epsilon}\int_{0}^{1}\rho^{2\gamma-(\theta+1)}dx\\
&&+\frac{1}{2}\int_{0}^{1}\rho^{\theta+1}\left(\partial_{x}(r^{n-1}u_0)\right)^{2}dx
   +\frac{1}{2}\int_{0}^{1}\rho^{2\gamma-(\theta+1)}dx
   +\int_{0}^{1}|(u-u_{0})f|dx.
\end{eqnarray*}
Choosing $\epsilon=\min\{c_2+\frac{2}{n}c_1,
\frac{2(n-1)}{n}c_1\}$, using $n\geq 2, (\ref{upalpha}),
(\ref{rhobdry})$ and $\rho_0^{1+\theta}(u_{0})_x^2\in L^1$, we get
\begin{eqnarray*}
                &&\frac{d}{dt}\int_{0}^{1}\frac{1}{2}(u-u_{0})^{2}dx
                +\int_{0}^{1}(c_{2}+\frac{2}{n}c_1-\epsilon)\rho^{\theta+1}
    \left(\partial_{x}(r^{n-1}u)\right)^{2}dx\\
    &&+\int_{0}^{1}(\frac{2(n-1)}{n}c_{1}-\epsilon)
    \rho^{\theta+1}\left(r^{n-1}\partial_{x}u-\frac{u}{r\rho}\right)^{2}dx\\
            &\leq& C_{\epsilon}\int_{0}^{1}\rho^{\theta+1}
            \left(\partial_{x}(r^{n-1}u_{0})\right)^{2}dx
            +C_{\epsilon}\int_{0}^{1}\rho^{\theta-1}u_{0}^{2}r^{-2}dx\\
                &&+C_{\epsilon}\int_{0}^{1}\rho^{\theta+1}r^{2n-2}
                (u_{0})_x^{2}dx+C_{\epsilon}\int_{0}^{1}\rho^{2\gamma-(\theta+1)}dx
                +\int_{0}^{1}|(u-u_{0})f|dx\\
    &\leq& C.
\end{eqnarray*}
 Thus, we can get (\ref{2.16}) immediately.
\end{proof}
Now we turn to prove \textbf{Corollary \ref{cor}}:
\begin{proof}
Using integration by parts, we get
\begin{eqnarray*}
&&\left|\int_{0}^{t}\int_{x}^{1}\frac{\partial_{t}u}{r^{n-1}}dxds\right|
=\left|\int_{x}^{1}\frac{u}{r^{n-1}}dx-\int_{x}^{1}\frac{u_0{}}{r_{0}^{n-1}}dx
   +\int_{0}^{t}\int_{x}^{1}(n-1)u^{2}r^{-n}dxds\right|\\
&\leq&\int_{x}^{1}\frac{|u-u_{0}|}{r^{n-1}}dx
      +\int_{x}^{1}\left|u_{0}(\frac{1}{r^{n-1}}-\frac{1}{r_{0}^{n-1}})\right|dx
      +\int_{0}^{t}\int_{x}^{1}(n-1)u^{2}r^{-n}dxds.\\
\end{eqnarray*}
Then, using Lemmas \ref{2.2}$\sim$\ref{ul2} and H\"{o}lder's
inequality, we have
\begin{eqnarray*}
&&(n-1)\int_{0}^{t}\int_{x}^{1}u^{2}r^{-n}dxds
\leq(n-1)a^{-n}\int_{0}^{t}\int_{x}^{1}u^{2}dxds\\
&\leq&(n-1)a^{-n}\int_{0}^{t}(\int_{0}^{1}u^{4m}dx)^{\frac{1}{2m}}dt(1-x)^{\frac{2m-1}{2m}}
\leq Ct(1-x)^{\frac{2m-1}{2m}},
\end{eqnarray*}
\begin{eqnarray*}
&&\int_{x}^{1}\frac{|u(y, t)-u_{0}(y)|}{r^{n-1}(y, t)}dy
\leq a^{1-n}(\int_{0}^{1}|u(y, t)-u_{0}(y)|^{2m}dy)^{\frac{1}{2m}}(1-x)^{\frac{2m-1}{2m}}\\
&\leq&a^{1-n}(\int_{0}^{1}|u(y,
t)-u_{0}(y)|^{2}dy)^{\frac{1}{2(2m-1)}}\\
&&\times
(\int_{0}^{1}|u(y, t)-u_{0}(y)|^{4m}dy)^{\frac{m-1}{2m(2m-1)}}(1-x)^{\frac{2m-1}{2m}}\\
&\leq&Ct^{\frac{1}{2(2m-1)}}(1-x)^{\frac{2m-1}{2m}},
\end{eqnarray*}
and
\begin{eqnarray*}
&&\int_{x}^{1}\left|u_{0}(\frac{1}{r^{n-1}}-\frac{1}{r_{0}^{n-1}})\right|dx
\leq\|u_{0}\|_{L^{\infty}_{x}}a^{2-2n}\int_{x}^{1}|r^{n-1}-r_{0}^{n-1}|dx\\
&\leq&\|u_{0}\|_{L^{\infty}_{x}}a^{2-2n}\int_{x}^{1}\left|\int_{0}^{t}\partial_tr^{n-1}ds\right|dx\\
&\leq& (n-1)\|u_{0}\|_{L^{\infty}_{x}}a^{2-2n}R^{n-2}
       \int_{x}^{1}\int_{0}^{t}|u|dsdx\\
&\leq&C\|u_{0}\|_{L^{\infty}_{x}}a^{2-2n}R^{n-2}
       \int_{0}^{t}(\int_{0}^{1}u^{4m}dx)^{\frac{1}{4m}}(1-x)^{\frac{4m-1}{4m}}ds
\leq Ct(1-x)^{\frac{4m-1}{4m}}.
\end{eqnarray*}
Finally, we can get (\ref{2.21}) immediately.
\end{proof}
\subsection{Proof of Lemma \ref{weightedrho}}
\begin{proof}Multiplying (\ref{2.3}) by $(1-x)^{\alpha_{0}}r^{\beta}\partial_{x}\rho^{\theta}$
  and integrating  it over $x$ from 0 to 1, using Young's inequality, we have
\begin{eqnarray*}
&&\int_{0}^{1}(1-x)^{\alpha_{0}}(r^{\beta}\partial_{x}\rho^{\theta})^{2}dx\\
&\leq&\frac{1}{2}
 \int_{0}^{1}(1-x)^{\alpha_{0}}(r^{\beta}\partial_{x}\rho^{\theta})^{2}dx
+C\int_{0}^{1}
  (1-x)^{\alpha_{0}}
  (r_{0}^{\beta}(\rho_{0}^{\theta})_x)^{2}dx\\
&&+C\int_{0}^{1}(1-x)^{\alpha_{0}}(r^{1+\beta-n}u)^{2}dx
+C\int_{0}^{1}(1-x)^{\alpha_{0}}(r_{0}^{1+\beta-n}u_{0})^{2}dx\\
&&+C\int_{0}^{1}
  (1-x)^{\alpha_{0}}
  \left(
   \int_{0}^{t}
       (
        -r^{\beta}\partial_{x}\rho^{\gamma}
        +r^{1+\beta-n}f
        +(1+\beta-n)u^{2}r^{\beta-n}
        )ds
   \right)^{2}dx.
\end{eqnarray*}
From Lemma \ref{u4m}, $\alpha_0>1-2\theta\alpha$ and
$m>\frac{1}{2(1-\theta\alpha)}$, we have
\begin{eqnarray*}
\int_0^1(1-x)^{\alpha_0}(r^{1+\beta-n}u)^2dx\leq
C\int_0^1(1-x)^{\frac{2m}{2m-1}\alpha_0}dx+C\int_0^1u^{4m}dx\leq
C,
\end{eqnarray*}
and
\begin{eqnarray*}
&&\int_0^1(1-x)^{\alpha_0}\left(\int_0^tu^2r^{\beta-n}(1+\beta-n)ds\right)^2dx\\
&\leq& C\int_0^1(1-x)^{\alpha_0}(\int_0^tu^2ds)^2dx \leq
C\int_0^1(1-x)^{\alpha_0}\int_0^tu^4dsdx\\
&\leq&C\int_0^1(1-x)^{\frac{m}{m-1}\alpha_0}dx+C\int_0^1\int_0^tu^{4m}dsdx
\leq C.
\end{eqnarray*}
Thus, we get
\begin{eqnarray*}
&&\int_{0}^{1}(1-x)^{\alpha_{0}}(r^{\beta}\partial_{x}\rho^{\theta})^{2}dx
\leq
C+C\int_{0}^{1}(1-x)^{\alpha_{0}}\int_{0}^{t}|r^{\beta}\partial_{x}\rho^{\gamma}|^{2}(x, s)dsdx\\
&\leq&
C+C\int_{0}^{t}\max(\rho^{\gamma-\theta})^{2}\int_{0}^{1}(1-x)^{\alpha_{0}}(r^{\beta}\partial_{x}\rho^{\theta})^{2}dxds\\
&\leq&
C+C\int_{0}^{t}\int_{0}^{1}(1-x)^{\alpha_{0}}(r^{\beta}\partial_{x}\rho^{\theta})^{2}dxds.
\end{eqnarray*}
Using Gronwall's inequality,  we can obtain (\ref{local-E3.13})
immediately..
\end{proof}
\subsection{Proof of lemma \ref{ut2}}
At first, we consider the following two lemmas.
\begin{lem}\label{ss-L3.11}
Under the assumptions of  Lemma \ref{u4m}, we have, for any
$t\in[0, T']$,
            \begin{eqnarray}
             &&\int^1_0\left(\rho^{\theta-1}
              u^2+\rho^{\theta+1}(\partial_xu)^2
              \right)dx+\int^t_0\int^1_0(\partial_tu)^2(x,s)dxds\nonumber\\&\leq&
              C+\int^t_0\int^1_0\rho^{3+\theta}(\partial_xu)^4dxds.\label{ss-E2.17-9}
            \end{eqnarray}
\end{lem}
\begin{proof}
Multiplying (\ref{1.9})$_2$ by $\partial_tu$, integrating it over
$[0,1]\times[0,t]$, and using the boundary conditions
(\ref{1.11})-(\ref{1.12}), we have
    \begin{eqnarray}
      &&\int^t_0\int^1_0(\partial_tu)^2(x,s)dxds\nonumber\\
            &=&\int^t_0\int^1_0\rho^\gamma\partial_x(r^{n-1}(\partial_tu))dxds\nonumber\\
                &&-\int^t_0\int^1_0
      (2c_1+c_2)\rho^{1+\theta}\partial_x(r^{n-1}u)
      \partial_x(r^{n-1}(\partial_tu))dxds\nonumber\\
      &&+\int^t_0\int^1_02c_1(n-1)\rho^\theta\partial_x(r^{n-2}u(\partial_tu))
      dxds-\int^t_0\int^1_0f(\partial_tu)dxds\nonumber\\
            &:=&\sum^4_{i=1}H_i.\label{ss-E2.18-9}
    \end{eqnarray}
Using (A3), (\ref{2.1}), (\ref{rhobdry}), (\ref{rup}),
Cauchy-Schwarz inequality, Young inequality and
$m>\frac{3}{4(1+\alpha\theta-\alpha)}$, we have
    \begin{eqnarray}
      &&H_2+H_3\nonumber\\
            &=&\left.\left\{\int^1_0\left[-\frac{2c_1+c_2}{2}\rho^{1+\theta}[\partial_x(r^{n-1}u)]^2
            +c_1(n-1)\rho^\theta\partial_x(r^{n-2}u^2)\right]dx\right\}\right|^t_0\nonumber\\
            &&+\int^t_0\int^1_0\{
                    (2c_1+c_2)(n-1)\rho^{1+\theta}\partial_x(r^{n-1}u)
                    \partial_x(r^{n-2}u^2)\nonumber\\
                    &&-
      \frac{(2c_1+c_2)}{2}(1+\theta)\rho^{2+\theta}[\partial_x(r^{n-1}u)]^3
      +2\theta
                c_1(n-1)\rho^{\theta+1}\frac{u}{r}[\partial_x(r^{n-1}u)]^2
                \nonumber\\
    &&-\theta
    c_1n(n-1)\rho^\theta\frac{u^2}{r^2}\partial_x(r^{n-1}u)+
            2nc_1(n-1)(n-2)\rho^{\theta-1}\frac{u^3}{r^3}\nonumber\\
    &&-3c_1(n-1)(n-2)\rho^\theta
            \frac{u^2}{r^2}\partial_x(r^{n-1}u)\}dxds\nonumber\\
                        &\leq&C+C\int^t_0\int^1_0\left[
                        \rho^{2+\theta}(\partial_xu)^3+\rho^{\theta-1}u^3\right]dxds\nonumber\\
                        &&-C^{-1}\int^1_0\left[\rho^{1+\theta}(\partial_xu)^2+\rho^{\theta-1}u^2\right]dx\nonumber\\
    &\leq&C+C\int^t_0\int^1_0\left[u^{4m}+\rho^{\frac{4m(\theta-1)}{4m-3}}+\rho^{3+\theta}(\partial_xu)^4
                        \right]dxds\nonumber\\
                        &&-C^{-1}\int^1_0\left[\rho^{1+\theta}(\partial_xu)^2+\rho^{\theta-1}u^2\right]dx\nonumber\\
    &\leq&C+C\int^t_0\int^1_0\rho^{3+\theta}(\partial_xu)^4
       dxds\nonumber\\
       &&-C^{-1}\int^1_0\left[\rho^{1+\theta}(\partial_xu)^2+\rho^{\theta-1}u^2\right]dx,
    \end{eqnarray}
        \begin{eqnarray}
          H_1
          &=&\left.\left\{\int^1_0\rho^\gamma\partial_x(r^{n-1}u)dx\right\}\right|^t_0
          +\int^t_0\int^1_0 \gamma\rho^{\gamma+1}[\partial_x(r^{n-1}u)]^2dxds\nonumber\\
                        &&-\int^t_0\int^1_02 (n-1)\rho^{\gamma}\frac{u}{r}\partial_x(r^{n-1}u)dxds
                        +\int^t_0\int^1_0 n(n-1)\rho^{\gamma-1}\frac{u^2}{r^2}dxds \nonumber\\
          &\leq&\int^1_0 \rho^\gamma\partial_x(r^{n-1}u)dx+C,
        \end{eqnarray}
        \begin{equation}
          H_4\leq\frac{1}{2}\int^t_0\int^1_0(\partial_tu)^2dxds+C.
          \qquad\qquad\qquad\qquad\qquad\qquad\label{ss-E2.22-9}
        \end{equation}
Using (\ref{ss-E2.18-9})-(\ref{ss-E2.22-9}), we can obtain
(\ref{ss-E2.17-9}) immediately.
\end{proof}
\begin{lem}\label{uxut}Under the assumptions of Lemma \ref{u4m},  there exist two constants
$C_{20}$ and $C_{21}$,  such that, for any $t\in[0, T']$
\begin{eqnarray}
\int_{0}^{1}\rho^{\theta+3}(\partial_{x}u)^{4}dx\leq
C_{20}+C_{21}\left(\int_{0}^{1}\left(\rho^{1+\theta}(\partial_xu)^2+
        \rho^{\theta-1}u^2+(\partial_tu)^{2}\right)dx\right)^2,\label{2.24}\\
\int_{0}^{1}(\rho\partial_{x}u)^{4}dx\leq
C_{20}+C_{21}\left(\int_{0}^{1}\left(\rho^{1+\theta}(\partial_xu)^2+
        \rho^{\theta-1}u^2+(\partial_tu)^{2}\right)dx\right)^{2}.\label{2.25}
  \end{eqnarray}
\end{lem}
\begin{proof}From (\ref{1.9})$_2$ and (\ref{1.12}), we have
    \begin{eqnarray}\label{**}
        \partial_{x}u
        &=&\frac{r^{1-n}}{2c_{1}+c_{2}}\rho^{\gamma-\theta-1}
        -(n-1)\rho^{-1}ur^{-n}
            -\frac{\rho^{-\theta-1}r^{1-n}}{2c_{1}+c_{2}}\int_{x}^{1}\frac{\partial_{t}u}{r^{n-1}}dy\nonumber\\
                    &&+\frac{2c_{1}(n-1)}{2c_{1}+c_{2}}\left(\frac{u}{\rho r^n}+\rho^{-\theta-1}
                    r^{1-n}\int_{x}^{1}\rho^{\theta}\partial_{y}(\frac{u}{r})dy
                    \right)\nonumber\\
        &&+\frac{\rho^{-\theta-1}r^{1-n}}{2c_{1}+c_{2}}\int_{x}^{1}\frac{f}{r^{n-1}}dy.
\end{eqnarray}
Using $r\geq a$, $\gamma>1$ and $\gamma>\theta$, we get
\begin{eqnarray*}
&&\int_0^1\rho^{\theta+3}(\partial_xu)^4dx\\
&\leq&C+C\int_0^1\rho^{\theta-1}dx +C\int_0^1\rho^{\theta-1}u^4dx
+C\int_0^1\rho^{-3\theta-1}\left(\int_x^1|\partial_tu|dy\right)^4dx\\
&&+C\int_0^1\!\rho^{-3\theta-1}\left(\int_x^1|\rho^{\theta}\partial_yu|\!+\!
|\rho^{\theta-1}u|dy\right)^4dx
\!+\!C\int_0^1\!\rho^{-3\theta-1}\left(\int_x^1|f|dy\right)^4dx.
\end{eqnarray*}
Then, using $\alpha(\theta-1)>-1$,
$\alpha<\min\{\frac{3}{3+\theta},\frac{3}{3\theta+1}\}$ and
$m>\frac{1}{1+\alpha\theta-\alpha}$, we have
$$
\int_0^1\rho^{\theta-1}dx\leq
C\int_0^1(1-x)^{\alpha(\theta-1)}dx\leq C,
$$
$$
\int_0^1\rho^{\theta-1}u^4dx \leq
C\int_0^1\rho^{(\theta-1)\frac{m}{m-1}}dx+C\int_0^1u^{4m}dx\leq C,
$$
\begin{eqnarray*}
&&\int_0^1\rho^{-3\theta-1}\left(\int_x^1|\partial_tu|dy\right)^4dx
\leq C\int_{0}^{1}\rho^{-3\theta-1}\left((\int_{0}^{1}|\partial_{t}u|^{2}dy)^{\frac{1}{2}}(1-x)^{\frac{1}{2}}\right)^{4}dx\\
&\leq&C\left(\int_{0}^{1}|\partial_{t}u|^{2}dx\right)^2\int_{0}^{1}(1-x)^{-\alpha(3\theta+1)+2}dx
\leq C\left(\int_{0}^{1}|\partial_{t}u|^{2}dx\right)^2,
\end{eqnarray*}
    $$
    \int_0^1\rho^{-3\theta-1}\left(\int_x^1|f|dy\right)^4dx \leq
    C\int_0^1(1-x)^{-\alpha(3\theta+1)+4}dx\leq C,
    $$
    \begin{eqnarray*}
        &&\int_0^1\rho^{-3\theta-1}\left(\int_x^1|\rho^{\theta}\partial_yu|dy\right)^4dx\\
                &\leq&C\left(\int^1_0\rho^{1+\theta}(\partial_xu)^2dx\right)^2\int^1_0(1-x)^{-\alpha(1+3\theta)}\left(
                \int^1_x(1-y)^{\alpha(\theta-1)}dy\right)^2dx\\
     &\leq& C\left(\int^1_0\rho^{1+\theta}(\partial_xu)^2dx\right)^2\int^1_0(1-x)^{2-\alpha(3+\theta)}dx\\
                &\leq& C\left(\int^1_0\rho^{1+\theta}(\partial_xu)^2dx\right)^2,
              \   \textrm{ when }\ \theta<1,
    \end{eqnarray*}
            \begin{eqnarray*}
              &&\int^1_0\rho^{-3\theta-1}\left(\int_x^1|\rho^{\theta-1}u|dy\right)^4dx\\
                    &\leq&C\left(\int^1_0\rho^{\theta-1}u^2dx\right)^2\int^1_0(1-x)^{-\alpha(1+3\theta)}\left(
                \int^1_x(1-y)^{\alpha(\theta-1)}dy\right)^2dx\\
     &\leq& C\left(\int^1_0\rho^{\theta-1}u^2dx\right)^2\int^1_0(1-x)^{2-\alpha(3+\theta)}dx\\
                &\leq& C\left(\int^1_0\rho^{\theta-1}u^2dx\right)^2,
             \   \textrm{ when }\ \theta<1,
            \end{eqnarray*}
and
    \begin{eqnarray*}
    &&\int_0^1\rho^{-3\theta-1}\left(\int_x^1|\rho^{\theta}\partial_yu|+|\rho^{\theta-1}u|dy\right)^4dx\\
            &\leq&
            C\left(\int^1_0\rho^{1+\theta}(\partial_xu)^2dx\right)^2+C\left(\int^1_0\rho^{\theta-1}u^2dx\right)^2,
            \   \textrm{ when }\ \theta\geq1.
    \end{eqnarray*}
Thus, we get (\ref{2.24}). Similarly, using
$\alpha<\frac{3}{4\theta}$ and $\alpha<\frac{3}{2(1+\theta)}$, we
can obtain (\ref{2.25}) easily.
\end{proof}

Now, we can prove \textbf{Lemma \ref{ut2}} as follows:
\begin{proof}
Differentiating $(\ref{1.9})_{2}$ respect to $t$,  we get
\begin{eqnarray}\label{2.27}
&&\partial_{tt}u\nonumber\\
    &=&\partial_{t}
         \left(           r^{n-1}\left(
                                 (2c_{1}+c_{2})\rho^{\theta+1}(r^{n-1}u)_x
                                -\rho^{\gamma}
                              \right)_x         \right)
         -2c_{1}(n-1)\partial_{t}\left(r^{n-2}u\partial_{x}\rho^{\theta}\right)+\partial_{t}f\nonumber\\
            & :=&I+J,
\end{eqnarray}
where
$I=r^{n-1}\partial_{x}\left((2c_{1}+c_{2})\rho^{\theta+1}\partial_{x}(r^{n-1}\partial_{t}u)\right)
  -2c_{1}(n-1)\left(r^{n-2}\partial_{t}u\partial_{x}\rho^{\theta}\right)$.

Multiplying $(\ref{2.27})$ by $\partial_{t}u$ and integrating it
over $x$ from 0 to 1, we have
\begin{eqnarray}\label{2.28}
 \frac{d}{dt}\int_{0}^{1}\frac{1}{2}(\partial_{t}u)^{2}dx
            +(\frac{2}{n}c_{1}+c_{2})\int_{0}^{1}\rho^{\theta+1}
               \left(\partial_{x}(r^{n-1}\partial_{t}u)\right)^2
              dx&&\nonumber\\
            +2\frac{(n-1)}{n}c_{1}\int_0^1\rho^{\theta+1}
                \left(r^{n-1}\partial_{tx}u-\rho^{-1}r^{-1}\partial_tu\right)^2dx
              &=&\int_{0}^{1}J\partial_{t}udx.\nonumber
\end{eqnarray}
Denote
\begin{eqnarray*}
B_{1}^{2}=\rho^{\theta+1}r^{2n-2}(\partial_{tx}u)^{2},\quad
 B_{2}^{2}=\rho^{\theta-1}r^{-2}(\partial_{t}u)^{2}.
\end{eqnarray*}
Thus
\begin{eqnarray*}
\int_0^1J\partial_tudx&\leq&\eta\int_0^1(B_1^2+B_2^2)dx
+C_{\eta}\int_0^1
         \left(
               \rho^{\theta+3}(\partial_xu)^4r^{4n-4}\right.\\
               &&+\rho^{\theta+1}(\partial_xu)^2u^2r^{2n-4}
             +\rho^{2\gamma-\theta+1}(\partial_xu)^2r^{2n-2}+\rho^2(\partial_xu)^2\\
           &&+\rho^{3-\theta}(\partial_xu)^2r^{2n-2}
           \left.+\rho^{\theta-1}u^4
         \right)dx
         +\int_0^1|\partial_tu\partial_tf|dx.
\end{eqnarray*}
Using $a\leq r\leq R$, $m>\frac{1}{1+\alpha\theta-\alpha}$,
$2\alpha(1-\theta)>-1$, Lemmas \ref{u4m} and \ref{uxut}, we get
$$
\int_0^1\rho^{\theta+3}(\partial_xu)^4r^{4n-4}dx \leq
C\left(\int_0^1\left(\rho^{1+\theta}(\partial_xu)^2+
        \rho^{\theta-1}u^2+(\partial_tu)^{2}\right)dx\right)^{2}+C,
$$
\begin{eqnarray*}
&&\int_0^1\rho^{\theta+1}(\partial_xu)^2u^2r^{2n-4}dx
\leq C\int_0^1\rho^{\theta+3}(\partial_xu)^4dx+C\int_0^1\rho^{\theta-1}u^4dx\\
&\leq& C\int_0^1\rho^{\theta+3}(\partial_xu)^4dx+C\int_0^1\rho^{(\theta-1)\frac{m}{m-1}}dx+C\int_0^1u^{4m}dx\\
&\leq& C\left(\int_0^1\left(\rho^{1+\theta}(\partial_xu)^2+
        \rho^{\theta-1}u^2+(\partial_tu)^{2}\right)dx\right)^{2}+C,
\end{eqnarray*}
\begin{eqnarray*}
&&\int_0^1\rho^{2\gamma-\theta+1}(\partial_xu)^2r^{2n-2}dx \leq
C\int_0^1\rho^{\theta+3}(\partial_xu)^4dx+C\int_0^1\rho^{4\gamma-3\theta-1}dx\\
&\leq& C\left(\int_0^1\left(\rho^{1+\theta}(\partial_xu)^2+
        \rho^{\theta-1}u^2+(\partial_tu)^{2}\right)dx\right)^{2}+C,
\end{eqnarray*}
\begin{eqnarray*}
&&\int_0^1\rho^{3-\theta}(\partial_xu)^2r^{2n-2}dx\leq
C\int_0^1\rho^{4}(\partial_xu)^4dx+C\int_0^1\rho^{2-2\theta}dx\\
&\leq& C\left(\int_0^1\left(\rho^{1+\theta}(\partial_xu)^2+
        \rho^{\theta-1}u^2+(\partial_tu)^{2}\right)dx\right)^{2}+C,
\end{eqnarray*}
\begin{eqnarray*}
&&\int_0^1\rho^2(\partial_xu)^2dx \leq
C\int_0^1\rho^{\theta+3}(\partial_xu)^4dx+C\int_0^1\rho^{1-\theta}dx\\
&\leq& C\left(\int_0^1\left(\rho^{1+\theta}(\partial_xu)^2+
        \rho^{\theta-1}u^2+(\partial_tu)^{2}\right)dx\right)^{2}+C,
\end{eqnarray*}
and
$$
\int_0^1\rho^{\theta-1}u^4dx\leq
C\int_0^1\rho^{(\theta-1)\frac{m}{m-1}}dx+C\int_0^1u^{4m}dx\leq C.
$$
Since
$\int_{0}^{1}\{r^{n-1}\partial_{x}[(2c_{1}+c_{2})\rho^{\theta+1}\partial_{x}(r^{n-1}u)-\rho^{\gamma}]
              -2c_{1}(n-1)r^{n-2}u\partial_{x}\rho^{\theta}+f\}^{2}dx|_{t=0}$
is bounded (see (A3)), we have $\int_0^1u_t^2(x,0)dx$ is bounded.
Thus, we have
\begin{eqnarray*}
&&\int_{0}^{1}(\partial_{t}u)^{2}dx+\int_0^t\int_0^1\left[\rho^{\theta+1}r^{2n-2}(\partial_{tx}u)^2
+\rho^{\theta-1}r^{-2}(\partial_tu)^2\right]dxds\\
&&\leq
C+C\int_0^t\left(\int_{0}^{1}\left(\rho^{1+\theta}(\partial_xu)^2+
        \rho^{\theta-1}u^2+(\partial_tu)^{2}\right)dx\right)^{2}ds.
\end{eqnarray*}
From (\ref{ss-E2.17-9}) and (\ref{2.24}), we have
\begin{eqnarray*}
&&\int_{0}^{1}\left(\rho^{1+\theta}(\partial_xu)^2+
        \rho^{\theta-1}u^2+(\partial_tu)^{2}\right)dx\\
&&+\int_0^t\int_0^1\left[\rho^{\theta+1}r^{2n-2}(\partial_{tx}u)^2+\rho^{\theta-1}r^{-2}(\partial_tu)^2\right]dxds\\
&\leq&
C_{24}+C_{25}\int_0^t\left(\int_{0}^{1}\left(\rho^{1+\theta}(\partial_xu)^2+
        \rho^{\theta-1}u^2+(\partial_tu)^{2}\right)dx\right)^{2}ds.
\end{eqnarray*}
Using Gronwall's inequality, we get
\begin{eqnarray}\label{C12}
&&\int_{0}^{1}\left(\rho^{1+\theta}(\partial_xu)^2+
        \rho^{\theta-1}u^2+(\partial_tu)^{2}\right)dx\nonumber\\
&&+\int_0^t\int_0^1\left[\rho^{\theta+1}r^{2n-2}(\partial_{tx}u)^2
+\rho^{\theta-1}r^{-2}(\partial_tu)^2\right]dxds\nonumber\\
&\leq&\frac{1}{\frac{1}{C_{24}}-C_{25}t} \leq 2C_{24}:=C_{6},
\indent0\leq t\leq\overline{T}_{4},
\end{eqnarray}
where
\begin{equation}\label{t1}
\overline{T}_{4}=\frac{1}{2C_{24}C_{25}}.
\end{equation}
Using Lemmas \ref{ss-L3.11}$\sim$\ref{uxut} and (\ref{C12}), we
can get (\ref{2.155})$\sim$(\ref{2.15}) immediately.
\end{proof}
\subsection{Proof of Lemma \ref{uxl1}}
\begin{proof}From (\ref{**}), we have
\begin{eqnarray*}
&&\left|\partial_{x}u\right|^{\lambda_0} \leq
C\left(\rho^{\lambda_0(\gamma-\theta-1)}
          +\left|\rho^{-1}u\right|^{\lambda_0}\right.
+
  \left|
  \rho^{-\theta-1}\int_{x}^{1}|\partial_{t}u|dx
  \right|
  ^{\lambda_0}\\
&&+\left|
 \rho^{-\theta-1}\int_{x}^{1}|f|dx
 \right|
 ^{\lambda_0}
\left.+\left|
     \rho^{-\theta-1}\int_{x}^{1}|\rho^{\theta}\partial_{x}u|dx
   \right|
   ^{\lambda_0}+\left|\rho^{-\theta-1}\int_x^1|u\rho^{\theta-1}|dx\right|^{\lambda_0}\right),\\
\end{eqnarray*}
for some positive constant $C$.

From (\ref{upalpha}), (\ref{rhobdry}), (\ref{zeta}), Lemmas
\ref{u4m} and \ref{ut2} we have
\begin{eqnarray*}
&&\int_{0}^{1}|\rho|^{\lambda_0(\gamma-\theta-1)}dx\leq
 C\int_{0}^{1}(1-x)^{\lambda_0\alpha(-\theta-1)}dx\leq
C,
\\&&
\int_{0}^{1}|\rho^{-1}u|^{\lambda_0}dx\leq\left(\int_{0}^{1}\rho^{-\frac{4m}{4m-\lambda_0}\lambda_0}dx\right)^{\frac{4m-\lambda_0}{4m}}
\left(\int_{0}^{1}u^{4m}dx\right)^{\frac{\lambda_0}{4m}}\leq C,
\\
&&\int_{0}^{1}\left|\rho^{-\theta-1}\int_{x}^{1}|\partial_{t}u
|dy\right|^{\lambda_0}dx
\leq\int_{0}^{1}|\rho^{-\theta-1}(1-x)^{\frac{1}{2}}|^{\lambda_0}dx\|\partial_{t}u\|^{\lambda_0}_{L^{2}}\\
&\leq&
C\int_{0}^{1}(1-x)^{(-\alpha(\theta+1)+\frac{1}{2})\lambda_0}dx\leq
C,
\\&&
 \!\!\!\!\!\!\!\!\!\!\!\!\!\!\!  \int_{0}^{1}\!\left|\rho^{-\theta-1}\int_{x}^{1}|\rho^{\theta}\partial_{x}u|dy\right|^{\lambda_0}dx
            \leq\int_{0}^{1}\!\left|\rho^{-\theta-1}
            (\int_{0}^{1}\!\rho^{\theta+3}(\partial_xu)^4dy)^{\frac{1}{4}}
            (\int_{x}^{1}\!\rho^{\theta-1}dy)^{\frac{3}{4}}\right|^{\lambda_0}\!\!dx\\
    &&\qquad\qquad\qquad\qquad\qquad\quad\leq C\int_{0}^{1}(1-x)^{-\lambda_0\alpha(\theta+1)}dx
    \leq C,
    \end{eqnarray*}
   and
       \begin{eqnarray*}
   \int_{0}^{1}\left|\rho^{-\theta-1}\int_{x}^{1}|\rho^{\theta-1}u|dy\right|^{\lambda_0}dx
            &\leq&\int_{0}^{1}\left|\rho^{-\theta-1}
            (\int_{0}^{1}\rho^{\theta-1}u^2dy)^{\frac{1}{2}}
            (\int_{x}^{1}\rho^{\theta-1}dy)^{\frac{1}{2}}\right|^{\lambda_0}dx\\
    &\leq&C\int_{0}^{1}(1-x)^{-\lambda_0\alpha(\theta+1)}dx
    \leq C,
    \end{eqnarray*}
where we use the fact that $\lambda_0<\min\{\frac{4m}{4\alpha
m+1},\frac{1}{\alpha(\theta+1)}\}$.

From above all, we can find a positive constant $C_{8}$ such that
$\int_{0}^{1}|\partial_{x}u|^{\lambda_0}dx\leq C_{8}$. Using Lemma
\ref{u4m}, (\ref{c15}) and Sobolev's embedding theorem, we can get
(\ref{b2}) immediately.
\end{proof}
\section*{Acknowledgment}
The authors would like to thank Professor Daoyuan Fang very much
for his guidance and helpful discussions.

\end{document}